\documentclass[11pt]{amsart}

\scrollmode
\usepackage{times}
\usepackage{amsmath,graphicx}
\usepackage{latexsym}
\usepackage{amscd,amsthm,amssymb}
\usepackage[all]{xypic}

\newtheorem{thm}{Theorem}
\newtheorem{prop}[thm]{Proposition}
\newtheorem{lem}[thm]{Lemma}
\newtheorem{cor}[thm]{Corollary}
\newtheorem{quest}[thm]{Question}

\theoremstyle{definition}
\newtheorem{ex}[thm]{Example}
\newtheorem{rem}[thm]{Remark}
\newtheorem{defn}[thm]{Definition}

\title[Geography of symplectic 4-manifolds with divisible canonical class]{On the geography of symplectic 4-manifolds with divisible canonical class}
\author{M.~J.~D.~Hamilton}
\email{mark.hamilton@math.lmu.de}
\date{January 19, 2009; MSC 2000: 57R17, 57N13, 14J29.}

\begin{document}

\begin{abstract} In this article we consider a version of the geography question for simply-connected symplectic 4-manifolds that takes into account the divisibility of the canonical class as an additional parameter. We also find new examples of 4-manifolds admitting several symplectic structures, inequivalent under deformation and self-diffeomorphisms of the manifold.
\end{abstract}

\maketitle

\tableofcontents

In this article we are interested in the geography of simply-connected closed symplectic 4-manifolds whose canonical classes have a given divisibility. The geography question in general aims at finding for any given pair of integers $(x,y)$ a closed 4-manifold $M$ with some {\it a priori} specified properties (e.g.~irreducible, spin, simply-connected, symplectic or complex) such that the Euler characteristic $e(M)$ equals $x$ and the signature $\sigma(M)$ equals $y$. This question has been considered for simply-connected symplectic 4-manifolds both in the spin and non-spin case for example in \cite{Go, DPZ, Pa1, Pa2}. Some further references are \cite{Ch, FSgeo, Pe, PPX}. We are interested in the same question for simply-connected symplectic 4-manifolds whose canonical class, considered as an element in second cohomology with integer coefficients, is divisible by a given integer $d>1$. Since the canonical class is characteristic, the first case $d=2$ corresponds to the general case of spin symplectic 4-manifolds. 

Geography questions are often formulated in terms of the invariants $c_1^2$ and $\chi_h$ instead of $e$ and $\sigma$, which for smooth closed 4-manifolds are defined by 
\begin{align*}
c_1^2(M)&=2e(M)+3\sigma(M)\\
\chi_h(M)&=\tfrac{1}{4}(e(M)+\sigma(M)).
\end{align*}
For complex 4-manifolds these numbers have the same value as the square of the first Chern class and the holomorphic Euler characteristic, making the definitions consistent. 

The constructions used in this article depend on generalized fibre sums of symplectic manifolds, also known as Gompf sums or normal connected sums \cite{Go, McW}, in particular in the form of knot surgery \cite{FSknot} and a generalized version of knot surgery along embedded surfaces of higher genus \cite{FSfam}. Some details on the generalized fibre sum can be found in Section \ref{sect def gen fibre sum}.

In Sections \ref{sect homotop ellipt surf}, \ref{sect geo spin neg sigma c pos} and \ref{sect geo non-spin neg sigma c pos} we consider the case $c_1^2=0$ and the spin and non-spin cases for $c_1^2>0$ and negative signature, while the case $c_1^2<0$ is covered at the end of Section \ref{sect general restrict canonic div}. We did not try to consider the case of non-negative signature, since even without a restriction on the divisibility of the canonical class such simply-connected symplectic 4-manifolds are known to be difficult to find. 

As a consequence of these geography results there often exist at the same lattice point in the $(\chi_h,c_1^2)$--plane several simply-connected symplectic 4-manifolds whose canonical classes have pairwise different divisibilities. It is natural to ask whether the same smooth 4-manifold can admit several symplectic structures with canonical classes of different divisibilities. This question is considered in Sections \ref{sect inequiv symp str} and \ref{sect ex inequiv symp str}. The symplectic structures with this property are inequivalent under deformations and orientation preserving self-diffeomorphisms of the manifold.  Similar examples have been found before on homotopy elliptic surfaces by McMullen and Taubes \cite{McT}, Smith \cite{Smi} and Vidussi \cite{V}. Another application of the geography question to the existence of inequivalent contact structures on certain 5-manifolds can be found in \cite{MHthesis}.

In the final part of this article an independent construction of simply-connected symplectic 4-manifolds with divisible canonical class is given by finding complex surfaces of general type with divisible canonical class. The construction uses branched coverings over smooth curves in pluricanonical linear systems $|nK|$.

\subsection*{Acknowledgements} The material in this article is part of the author's Ph.D.~thesis, submitted in May 2008 at the University of Munich. I would like to thank D.~Kotschick, who supervised the thesis, as well as the {\it Studienstiftung} and the {\it DFG} for financial support.

\section{General restrictions on the divisibility of the canonical class}\label{sect general restrict canonic div}

We begin by deriving a few general restrictions for symplectic 4-manifolds admitting a symplectic structure whose canonical class is divisible by an integer $d>1$. 

Let $(M,\omega)$ be a closed, symplectic 4-manifold. The canonical class of the symplectic form $\omega$, denoted by $K$, is defined as
\begin{equation*}
K=-c_1(TM,J),
\end{equation*}
where $J$ is an almost complex structure compatible with $\omega$. The self-intersection number of $K$ is given by the formula
\begin{equation*}
K^2=c_1^2(M)=2e(M)+3\sigma(M).
\end{equation*}
Since the first Chern class $c_1(M,J)$ is characteristic, it follows by a general property of the intersection form that
\begin{equation*}
\mbox{$c_1^2(M)\equiv \sigma(M)$ mod $8$}
\end{equation*}
and hence the number 
\begin{equation*}
\chi_h(M)=\tfrac{1}{4}(e(M)+\sigma(M))
\end{equation*}
is an integer. If $b_1(M)=0$, this number is equal to $\frac{1}{2}(1+b_2^+(M))$. In particular, in this case $b_2^+(M)$ has to be an odd integer and $\chi_h(M)>0$. There is a further constraint if the manifold $M$ is spin, equivalent to the congruence $\sigma(M)\equiv 0$ mod $16$ given by Rochlin's theorem \cite{Ro1}:
\begin{equation*}
\mbox{$c_1^2(M)\equiv 8\chi_h(M)$ mod $16$}.
\end{equation*}
In particular, $c_1^2(M)$ is divisible by $8$. We say that $K$ is divisible by an integer $d$ if there exists a cohomology class $A\in H^2(M;\mathbb{Z})$ with $K=dA$. 

\begin{lem}\label{lem 1 geography K} Let $(M,\omega)$ be a closed symplectic 4-manifold. Suppose that the canonical class $K$ is divisible by an integer $d$. Then $c_1^2(M)$ is divisible by $d^2$ if $d$ is odd and by $2d^2$ if $d$ is even.
\end{lem}
\begin{proof} If $K$ is divisible by $d$, we can write $K=dA$, where $A\in H^2(M;\mathbb{Z})$. The equation $c_1^2(M)=K^2=d^2A^2$ implies that $c_1^2(M)$ is divisible by $d^2$ in any case. If $d$ is even, then $w_2(M)\equiv K\equiv 0$ mod $2$, hence $M$ is spin and the intersection form $Q_M$ is even. This implies that $A^2$ is divisible by $2$, hence $c_1^2(M)$ is divisible by $2d^2$. 
\end{proof}
The case $c_1^2(M)=0$ is special, since there are no restrictions from this lemma (see Section \ref{sect homotop ellipt surf}). For the general case of spin symplectic 4-manifolds ($d=2$) we recover the constraint that $c_1^2$ is divisible by $8$.

Further restrictions come from the adjunction formula
\begin{equation*}
2g-2=K\cdot C+C\cdot C,
\end{equation*}
where $C$ is an embedded symplectic surface of genus $g$, oriented by the restriction of the symplectic form. 
\begin{lem}\label{lem 2 geography K} Let $(M,\omega)$ be a closed symplectic 4-manifold. Suppose that the canonical class $K$ is divisible by an integer $d$.
\begin{itemize}
\item If $M$ contains a symplectic surface of genus $g$ and self-intersection $0$, then $d$ divides $2g-2$.
\item If $d\neq 1$, then $M$ is minimal. If the manifold $M$ is in addition simply-connected, then it is irreducible. 
\end{itemize}
\end{lem}
\begin{proof} The first part follows immediately by the adjunction formula. If $M$ is not minimal, then it contains a symplectically embedded sphere $S$ of self-intersection $(-1)$. The adjunction formula can be applied and yields $K\cdot S=-1$, hence $K$ has to be indivisible. The claim about irreducibility follows from \cite[Corollary 1.4]{MHDK}.
\end{proof}

The canonical class of a 4-manifold $M$ with $b_2^+\geq 2$ is a Seiberg-Witten basic class, i.e.~it has non-vanishing Seiberg-Witten invariant. This implies that only finitely many classes in $H^2(M;\mathbb{Z})$ can occur as the canonical classes of symplectic structures on $M$. The following is proved in \cite{LL3}.
\begin{thm}\label{canon classes equal b2+=1} Let $M$ be a (smoothly) minimal closed 4-manifold with $b_2^+=1$. Then the canonical classes of all symplectic structures on $M$ are equal up to sign.
\end{thm}

If $M$ is a K\"ahler surface, we can consider the canonical class of the K\"ahler form. 
\begin{thm}\label{thm geog intro Kahler} Suppose that $M$ is a minimal K\"ahler surface with $b_2^+>1$. 
\begin{itemize}
\item If $M$ is of general type, then $\pm K_M$ are the only Seiberg-Witten basic classes of $M$. 
\item If $N$ is another minimal K\"ahler surface with $b_2^+>1$ and $\phi\colon M\rightarrow N$ a diffeomorphism, then $\phi^*K_N=\pm K_M$.   
\end{itemize}
\end{thm}
For the proofs see \cite{FrM}, \cite{Mor} and \cite{Wi}. Choosing $\phi$ as the identity diffeomorphism, an immediate consequence of the second part of this theorem is:
\begin{cor} Let $M$ be a (smoothly) minimal closed 4-manifold with $b_2^+>1$. Then the canonical classes of all K\"ahler structures on $M$ are equal up to sign.
\end{cor}
The corresponding statement is not true in general for the canonical classes of {\em symplectic structures} on minimal 4-manifolds with $b_2^+>1$. There exist such 4-manifolds $M$ which admit several symplectic structures whose canonical classes in $H^2(M;\mathbb{Z})$ are not equal up to sign. In addition, such examples can be constructed where the canonical classes cannot be permuted by orientation-preserving self-diffeomorphisms of the manifold \cite{McT, Smi, V}, for example because they have different divisibilities as elements in integral cohomology (cf.~the examples in Sections \ref{sect inequiv symp str} and \ref{sect ex inequiv symp str}). 

It is useful to define the (maximal) divisibility of the canonical class in the case that $H^2(M;\mathbb{Z})$ is torsion free.  

\begin{defn}\label{defn divisib for abel group} Suppose $H$ is a finitely generated free abelian group. For $a\in H$ let
\begin{equation*}
d(a)=\max\{k\in\mathbb{N}_0\mid \mbox{there exists an element $b\in H$, $b\neq 0$, with $a=kb$}\}.
\end{equation*}
We call $d(a)$ the {\em divisibility} of $a$ (or sometimes, to emphasize, the {\em maximal} divisibility). With this definition the divisibility of $a$ is $0$ if and only $a=0$. We call $a$ {\em indivisible} if $d(a)=1$.  
\end{defn}

In particular, if $M$ is a simply-connected manifold the integral cohomology group $H^2(M;\mathbb{Z})$ is torsion free and the divisibility of the canonical class $K\in H^2(M;\mathbb{Z})$ is well-defined.

\begin{prop}\label{prop no complex surface has two symp with diff div} Suppose that $M$ is a simply-connected closed 4-manifold which admits at least two symplectic structures whose canonical classes have different divisibilities. Then $M$ is not diffeomorphic to a complex surface.
\end{prop}
\begin{proof} The assumptions imply that $M$ has a symplectic structure whose canonical class has divisibility $\neq 1$. By Lemma \ref{lem 2 geography K} the manifold $M$ is (smoothly) minimal and by Theorem \ref{canon classes equal b2+=1} it has $b_2^+>1$. Suppose that $M$ is diffeomorphic to a complex surface. The Kodaira-Enriques classification implies that $M$ is diffeomorphic to an elliptic surface $E(n)_{p,q}$ with $n\geq 2$ and $p,q$ coprime or to a surface of general type. 

Consider the elliptic surfaces $E(n)_{p,q}$ for $n\geq 2$ and denote the class of a general fibre by $F$. The Seiberg-Witten basic classes of these 4-manifolds are known \cite{FSrational}. They consist of the set of classes of the form $kf$ where $f$ denotes the indivisible class $f=F/pq$ and $k$ is an integer
\begin{equation*}
\mbox{$k\equiv npq-p-q$ mod $2$},\quad |k|\leq npq-p-q.
\end{equation*}
Suppose that $\omega$ is a symplectic structure on $E(n)_{p,q}$ with canonical class $K$. By a theorem of Taubes \cite{Tmore, Bourbaki} the inequality 
\begin{equation*}
K\cdot [\omega]\geq |c\cdot[\omega]|
\end{equation*}
holds for any basic class $c$, with equality if and only if $K=\pm c$, and the number $K\cdot [\omega]$ is positive if $K$ is non-zero. It follows that the canonical class of any symplectic structure on $E(n)_{p,q}$ is given by $\pm(npq-p-q)f$, hence there is only one possible divisibility. This follows for surfaces of general type by the first part of Theorem \ref{thm geog intro Kahler}.
\end{proof} 

We now consider the geography question for manifolds with $c_1^2<0$. The following theorem is due to C.~H.~Taubes \cite{swg} in the case $b_2^+\geq 2$ and to A.-K.~Liu \cite{Liu} in the case $b_2^+=1$.
\begin{thm}\label{taubes liu minimal K^2} Let $M$ be a closed, symplectic 4-manifold. Suppose that $M$ is minimal. 
\begin{itemize}
\item If $b_2^+(M)\geq 2$, then $K^2\geq 0$.
\item If $b_2^+(M)=1$ and $K^2<0$, then $M$ is a ruled surface, i.e.~an $S^2$-bundle over a surface (of genus $\geq 2$).
\end{itemize}
\end{thm}

Since ruled surfaces over irrational curves are not simply-connected, any simply-connected, symplectic 4-manifold $M$ with $c_1^2(M)<0$ is not minimal. By Lemma \ref{lem 2 geography K} this implies that $K$ is indivisible, $d(K)=1$. 

Let $(\chi_h,c_1^2)=(n,-r)$ be a lattice point with $n,r\geq 1$ and $M$ a simply-connected symplectic 4-manifold with these invariants. Since $M$ is not minimal, we can successively blow down $r$ $(-1)$-spheres in $M$ to get a simply-connected symplectic 4-manifold $N$ with invariants $(\chi_h,c_1^2)=(n,0)$ such that there exists a diffeomorphism $M=N\#r\overline{\mathbb{C}P^2}$. 

Conversely, consider the manifold 
\begin{equation*}
M=E(n)\#r\overline{\mathbb{C}P^2}.
\end{equation*}
Then $M$ is a simply-connected symplectic 4-manifold with indivisible $K$. Since $\chi_h(E(n))=n$ and $c_1^2(E(n))=0$, this implies 
\begin{equation*}
(\chi_h(M),c_1^2(M))=(n,-r).
\end{equation*}
Hence the point $(n,-r)$ can be realized by a simply-connected symplectic 4-manifold.

\section{The generalized fibre sum}\label{sect def gen fibre sum}

In this section we recall the definition of the generalized fibre sum from \cite{Go, McW} and fix some notation, used in \cite{MHthesis}. Let $M$ and $N$ be closed oriented 4-manifolds which contain embedded oriented surfaces $\Sigma_M$ and $\Sigma_N$ of genus $g$ and self-intersection $0$. We choose trivializations of the form $\Sigma_g\times D^2$ for tubular neighbourhoods of the surfaces $\Sigma_M$ and $\Sigma_N$. The generalized fibre sum
\begin{equation*}
X=M\#_{\Sigma_M=\Sigma_N}N
\end{equation*}
is then formed by deleting the interior of the tubular neighbourhoods and gluing the resulting manifolds $M'$ and $N'$ along their boundaries $\Sigma_g\times S^1$, using a diffeomorphism which preserves the meridians to the surfaces, given by the $S^1$--fibres, and reverses the orientation on them. The closed oriented 4-manifold can depend on the choice of trivializations and gluing diffeomorphism. The trivializations of the tubular neighbourhoods also determine push-offs of the central surfaces $\Sigma_M$ and $\Sigma_N$ into the boundary. Under inclusion the push-offs determine surfaces $\Sigma_X$ and $\Sigma_X'$ of self-intersection $0$ in the 4-manifold $X$. In general, these surfaces do not represent the same homology class in $X$ but differ by a rim torus. However, if the gluing diffeomorphism is chosen such that it preserves also the $\Sigma_g$--fibres in $\Sigma_g\times S^1$, then the push-offs get identified to a well-defined surface $\Sigma_X$ in $X$.

Suppose that the surfaces $\Sigma_M$ and $\Sigma_N$ represent indivisible non-torsion classes in the homology of $M$ and $N$. We can then choose surfaces $B_M$ and $B_N$ in $M$ and $N$ which intersect $\Sigma_M$ and $\Sigma_N$ in a single positive transverse point. These surfaces with a disk removed can be assumed to bound the meridians to $\Sigma_M$ and $\Sigma_N$ in the manifolds $M'$ and $N'$, hence they sew together to give a surface $B_X$ in $X$.

The second cohomology of $M$ can be split into a direct sum
\begin{equation*}
H^2(M;\mathbb{Z})\cong P(M)\oplus\mathbb{Z}\Sigma_M\oplus \mathbb{Z}B_M,
\end{equation*}
where $P(M)$ denotes the orthogonal complement to the subgroup $\mathbb{Z}\Sigma_M\oplus \mathbb{Z}B_M$ in $H^2(M;\mathbb{Z})$ with respect to the intersection form $Q_M$. The restriction of the intersection form to the last two summands is given by
\begin{equation*}
\left(\begin{array}{cc} 0 &1 \\ 1& B_M^2 \\ \end{array}\right).
\end {equation*}
This form is unimodular, hence the restriction of the intersection form to $P(M)$ (modulo torsion) is unimodular as well. There exists a similar decomposition for the second cohomology of $N$.

\begin{thm}\label{thm formula second cohom fibre sum} Suppose that the integral cohomology of $M$, $N$ and $X$ is torsion-free and the surfaces $\Sigma_M$ and $\Sigma_N$ represent indivisible classes. If rim tori do not exist in the fibre sum $X=M\#_{\Sigma_M=\Sigma_N}N$, then the second cohomology of $X$ splits as a direct sum
\begin{equation*}
H^2(X;\mathbb{Z})\cong P(X)\oplus \mathbb{Z}\Sigma_X\oplus\mathbb{Z}B_X,
\end{equation*}
where 
\begin{equation*}
P(X)\cong P(M)\oplus P(N).
\end{equation*}
The restriction of the intersection form $Q_X$ to $P(X)$ is the direct sum of the restrictions of $Q_M$ and $Q_N$ and the restriction to $\mathbb{Z}\Sigma_X\oplus\mathbb{Z}B_X$ is of the form 
\begin{equation*}
\left(\begin{array}{cc} 0 &1 \\ 1& B_M^2+B_N^2 \\ \end{array}\right).
\end {equation*}
\end{thm}
A proof for this theorem can be found in \cite[Section V.3.5]{MHthesis}. It implies that there exist monomorphisms of abelian groups of both $H^2(M;\mathbb{Z})$ and $H^2(N;\mathbb{Z})$ into $H^2(X;\mathbb{Z})$ given by
\begin{equation}\label{eq monomorph M into X}
\begin{split}
\Sigma_M&\mapsto \Sigma_X\\
B_M&\mapsto B_X\\
Id:P(M)&\rightarrow P(M),
\end{split}
\end{equation}
and similarly for $N$. The monomorphisms do not preserve the intersection form if $B_M^2$ or $B_N^2$ differ from $B_X^2$. The following lemma is sometimes useful in checking the conditions for Theorem \ref{thm formula second cohom fibre sum} (the proof follows from Sections V.2 and V.3 in \cite{MHthesis}).
\begin{lem}\label{lem no rim tori X isom M} Let $X=M\#_{\Sigma_M=\Sigma_N}N$ be a generalized fibre sum along embedded surfaces of self-intersection $0$. Suppose that the map on integral first homology induced by one of the embeddings, say $\Sigma_N\rightarrow N$, is an isomorphism. Then rim tori do not exist in $X$. If in addition one of the surfaces represents an indivisible homology class, then $H_1(X;\mathbb{Z})\cong H_1(M;\mathbb{Z})$.     
\end{lem}

Suppose that $M$ and $N$ are symplectic 4-manifolds and $\Sigma_M$ and $\Sigma_N$ symplectically embedded. We orient both surfaces by the restriction of the symplectic forms. Then the generalized fibre sum $X$ also admits a symplectic structure. The canonical class $K_X$ can be calculated by the following formula:

\begin{thm}\label{thm formula can class gen fibre sum} Under the assumptions of Theorem \ref{thm formula second cohom fibre sum} and the embeddings of the cohomology of $M$ and $N$ into the cohomology of $X$ given by equation \eqref{eq monomorph M into X}, we have
\begin{equation*}
K_X=K_M+K_N-(2g-2)B_X+2\Sigma_X.
\end{equation*}
\end{thm}
A proof can be found in \cite[Section V.5]{MHthesis}. The formula for $g=1$ has been proved in \cite{Smi} and a related formula for arbitrary $g$ can be found in \cite{IP2}.

\section{The knot surgery construction}\label{sect knot surg}

The following construction due to Fintushel and Stern \cite{FSknot} is used frequently in the following sections. Let $K$ be a knot in $S^3$ and denote a tubular neighbourhood of $K$ by $\nu K\cong S^1\times D^2$. Let $m$ be a fibre of the circle bundle $\partial \nu K\rightarrow K$ and use an oriented Seifert surface for $K$ to define a section $l\colon K\rightarrow\partial \nu K$. The circles $m$ and $l$ are called the meridian and the longitude of $K$. Let $M_K$ be the closed 3-manifold obtained by $0$-Dehn surgery on $K$. The manifold $M_K$ is constructed in the following way: Consider $S^3\setminus \mbox{int}\,\nu K$ and let 
\begin{equation*}
f\colon \partial (S^1\times D^2)\rightarrow \partial (S^3\setminus \mbox{int}\,\nu K)
\end{equation*}
be a diffeomorphism which maps the circle $\partial D^2$ onto $l$. Then one defines 
\begin{equation*}
M_K=(S^3\setminus \mbox{int}\,\nu K)\cup_f (S^1\times D^2).
\end{equation*}
The manifold $M_K$ is determined by this construction uniquely up to diffeomorphism. One can show that it has the same integral homology as $S^2\times S^1$. The meridian $m$, which bounds the fibre in the normal bundle to $K$ in $S^3$, becomes non-zero in the homology of $M_K$ and defines a generator for $H_1(M_K;\mathbb{Z})$. The longitude $l$ is null-homotopic in $M_K$, since it bounds one of the $D^2$--fibres glued in. This disk fibre together with the Seifert surface of $K$ determine a closed, oriented surface $B_K$ in $M_K$ which intersects $m$ once and generates $H_2(M_K;\mathbb{Z})$. 

We consider the closed, oriented 4-manifold $M_K\times S^1$. It contains an embedded torus $T_K=m\times S^1$ of self-intersection $0$, which has a framing coming from a canonical framing of $m$. Let $X$ be an arbitrary closed, oriented 4-manifold, which contains an embedded torus $T_X$ of self-intersection $0$, representing an indivisible homology class. Then the result of {\em knot surgery} on $X$ is given by the generalized fibre sum 
\begin{equation*}
X_K=X\#_{T_X=T_K}(M_K\times S^1).
\end{equation*}
Here we have implicitly chosen a trivialization of the form $T^2\times D^2$ for the tubular neighbourhood of the torus $T_X$. We choose a gluing diffeomorphism which preserves both the $T^2$--factor and the $S^1$--factor on the boundaries of the tubular neighbourhoods and reverses orientation on the $S^1$--factor (the smooth 4-manifold $X_K$ might depend on the choice of the framing for $T_X$). The embedded torus of self-intersection $0$ in $X_K$, defined by identifying the push-offs, is denoted by $T_{X_K}$.  

The closed surface $B_K$ in the 3-manifold $M_K$ determines under inclusion a closed surface in the 4-manifold $M_K\times S^1$, denoted by the same symbol. It intersects the torus $T_K$ in a single transverse point. We also choose a surface $B_X$ in $X$ intersecting $T_X$ transversely and geometrically once. Both surfaces sew together to determine a surface $B_{X_K}$ in $X_K$ which intersects the torus $T_{X_K}$ in a single transverse point. 

We assume that the cohomology of $X$ is torsion free. From \cite{FSknot} it is known that there exists an isomorphism 
\begin{equation}\label{knot surgery decomp H^2}
H^2(X;\mathbb{Z})\cong H^2(X_K;\mathbb{Z})
\end{equation}
preserving intersection forms. In the notation of Section \ref{sect def gen fibre sum} this follows, because
\begin{equation*}
H^2(M_K\times S^1;\mathbb{Z})\cong \mathbb{Z}T_K\oplus \mathbb{Z}B_K,
\end{equation*}
hence $P(M_K\times S^1)=0$. In addition, the self-intersection number of $B_{X_K}$ is equal to the self-intersection number of $B_X$, because the class $B_K$ has zero self-intersection (it can be moved away in the $S^1$ direction). The claim then follows from Theorem \ref{thm formula second cohom fibre sum} and Lemma \ref{lem no rim tori X isom M}.

In particular, assume that both $X$ and $X'=X\setminus T_X$ are simply-connected. Since the fundamental group of $M_K\times S^1$ is normally generated by the image of the fundamental group of $T_K$ under inclusion, it follows that $X_K$ is again simply-connected, hence by Freedman's theorem \cite{Freedm} the manifolds $X$ and $X_K$ are homeomorphic. However, one can show with Seiberg-Witten theory that $X$ and $X_K$ are in many cases not diffeomorphic \cite{FSknot}.  

Suppose that $K$ is a fibred knot, i.e.~there exists a fibration
\begin{equation*}
\begin{CD}
S^3\setminus \mbox{int}\,\nu K @<<< \Sigma_{h}'\\
@VVV \\
S^1
\end{CD}
\end{equation*}
over the circle, where the fibres $\Sigma_h'$ are punctured surfaces of genus $h$ forming Seifert surfaces for $K$. Then $M_K$ is fibred by closed surfaces $B_K$ of genus $h$. This induces a fibre bundle
\begin{equation*}
\begin{CD}
M_K\times S^1 @<<< \Sigma_{h}\\
@VVV \\
T^2
\end{CD}
\end{equation*}
and the torus $T_K=m\times S^1$ is a section of this bundle. By a theorem of Thurston \cite{Thu} the manifold $M_K\times S^1$ admits a symplectic form such that $T_K$ and the fibres are symplectic. This construction can be used to do symplectic generalized fibre sums along $T_K$ if the manifold $X$ is symplectic and the torus $T_X$ symplectically embedded. The canonical class of $M_K\times S^1$ can be calculated by the adjunction formula, because the fibres $B_K$ and the torus $T_K$ are symplectic surfaces and form a basis of $H_2(M_K\times S^1;\mathbb{Z})$:
\begin{equation*}
K_{M_K\times S^1}=(2h-2)T_K.
\end{equation*}
The canonical class of the symplectic 4-manifold $X_K$ is then given by \cite{FSknot} 
\begin{equation}\label{Fint Stern knot K}
K_{X_K}=K_X+2hT_X,
\end{equation}
cf.~Theorem \ref{thm formula can class gen fibre sum}.

\section{Symplectic 4-manifolds with $c_1^2=0$}\label{sect homotop ellipt surf}

Having covered the case $c_1^2<0$ already in Section \ref{sect general restrict canonic div}, we now consider the case $c_1^2=0$.
\begin{defn} A closed, simply-connected 4-manifold $M$ is called a {\em homotopy elliptic surface} if $M$ is homeomorphic to a relatively minimal, simply-connected elliptic surface, i.e.~to a complex surface of the form $E(n)_{p,q}$ with $p,q$ coprime and $n\geq 1$.
\end{defn}

For details on the surfaces $E(n)_{p,q}$ see \cite[Section 3.3]{GS}. By definition, homotopy elliptic surfaces $M$ are simply-connected and have invariants
\begin{align*}
c_1^2(M)&=0\\
e(M)&= 12n\\
\sigma(M)&= -8n.
\end{align*}
The integer $n$ is equal to $\chi_h(M)$. In particular, symplectic homotopy elliptic surfaces have $K^2=0$. There exists the following converse.

\begin{lem}\label{hom ellipt} Let $M$ be a closed, simply-connected, symplectic 4-manifold with $K^2=0$. Then $M$ is a homotopy elliptic surface. 
\end{lem}
\begin{proof} Since $M$ is almost complex, the number $\chi_h(M)$ is an integer. The Noether formula
\begin{equation*}
\chi_h(M)=\tfrac{1}{12}(K^2+e(M))=\tfrac{1}{12}e(M)
\end{equation*}
implies that $e(M)$ is divisible by $12$, hence $e(M)=12k$ for some $k>0$. Together with the equation 
\begin{equation*}
0=K^2=2e(M)+3\sigma(M)
\end{equation*}
it follows that $\sigma(M)=-8k$. Suppose that $M$ is non-spin. If $k$ is odd, then $M$ has the same Euler characteristic, signature and type as $E(k)$. If $k$ is even, then $M$ has the same Euler characteristic, signature and type as the non-spin manifold $E(k)_2$. Since $M$ is simply-connected, $M$ is homeomorphic to the corresponding elliptic surface by Freedman's theorem \cite{Freedm}.

Suppose that $M$ is spin. Then the signature is divisible by $16$ due to Rochlin's theorem. Hence the integer $k$ above has to be even. Then $M$ has the same Euler characteristic, signature and type as the spin manifold $E(k)$. Again by Freedman's theorem the 4-manifold $M$ is homeomorphic to $E(k)$.    
\end{proof}

\begin{lem} Suppose that $M$ is a symplectic homotopy elliptic surface such that the divisibility of $K$ is even. Then $\chi_h(M)$ is even. 
\end{lem}
\begin{proof} The assumption implies that $M$ is spin. The Noether formula then shows that $\chi_h(M)$ is even, since $K^2=0$ and $\sigma(M)$ is divisible by $16$. 
\end{proof}

The next theorem shows that this is the only restriction on the divisibility of the canonical class $K$ for symplectic homotopy elliptic surfaces.

\begin{thm}[Homotopy elliptic surfaces]\label{elliptic div} Let $n$ and $d$ be positive integers. If $n$ is odd, assume in addition that $d$ is odd. Then there exists a symplectic homotopy elliptic surface $(M,\omega)$ with $\chi_h(M)=n$ whose canonical class $K$ has divisibility equal to $d$.
\end{thm}
Note that there is no constraint on $d$ if $n$ is even.
\begin{proof} If $n$ is equal to $1$ or $2$, the symplectic manifold can be realized as an elliptic surface. The canonical class of an elliptic surface $E(n)_{p,q}$ with $p,q$ coprime is given by
\begin{equation*}
K=(npq-p-q)f,
\end{equation*}
where $f$ is indivisible and $F=pqf$ denotes the class of a general fibre. For $n=1$ and $d$ odd we can take the surface $E(1)_{d+2,2}$, since 
\begin{equation*}
(d+2)2-(d+2)-2=d.
\end{equation*}
For $n=2$ and $d$ arbitrary we can take $E(2)_{d+1}=E(2)_{d+1,1}$, since 
\begin{equation*}
2(d+1)-(d+1)-1=d.
\end{equation*}
We now consider the case $n\geq 1$ in general and separate the proof into several cases. Suppose that $\mathbf{d=2k}$ and $\mathbf{n=2m}$ are both even with $k,m\geq 1$. Consider the elliptic surface $E(n)$. It contains a general fibre $F$ which is an embedded symplectic torus of self-intersection $0$. In addition, it contains a rim torus $R$ which arises from a decomposition of $E(n)$ as a fibre sum $E(n)=E(n-1)\#_F E(1)$, cf.~\cite{GoM} and Example \ref{n-1 Lagrangian triples in E(n)}. The rim torus $R$ has self-intersection $0$ and there exists a dual (Lagrangian) 2-sphere $S$ with intersection $RS=1$. We can assume that $R$ and $S$ are disjoint from the fibre $F$. The rim torus is in a natural way Lagrangian. By a perturbation of the symplectic form we can assume that it becomes symplectic. We give $R$ the orientation induced by the symplectic form. The proof consists in doing knot surgery along the fibre $F$ and the rim torus $R$.

Let $K_1$ be a fibred knot of genus $g_1=m(k-1)+1$. We do knot surgery along $F$ with the knot $K_1$ to get a new symplectic 4-manifold $M_1$. The elliptic fibration $E(n)\rightarrow \mathbb{C}P^1$ has a section showing that the meridian of $F$, which is the $S^1$--fibre of $\partial\nu F\rightarrow F$, bounds a disk in $E(n)\setminus \mbox{int}\,\nu F$. This implies that the complement of $F$ in $E(n)$ is simply-connected, hence the manifold $M_1$ is again simply-connected. By the knot surgery construction the manifold $M_1$ is homeomorphic to $E(n)$. The canonical class is given by formula \eqref{Fint Stern knot K}:
\begin{align*}
K_{M_1}&=(n-2)F+2g_1F\\
&=(2m-2+2mk-2m+2)F\\
&=2mkF.
\end{align*}
Here we have identified the cohomology of $M_1$ and $E(n)$ under the isomorphism in equation \eqref{knot surgery decomp H^2}. The rim torus $R$ still exists as an embedded oriented symplectic torus in $M_1$ with a dual 2-sphere $S$, because we can assume that the knot surgery takes place in a small neighbourhood of $F$ disjoint from $R$ and $S$. In particular, the complement of $R$ in $M_1$ is simply-connected. Let $K_2$ be a fibred knot of genus $g_2=k$ and $M$ the result of knot surgery on $M_1$ along $R$. Then $M$ is a simply-connected symplectic 4-manifold homeomorphic to $E(n)$. The canonical class of $M$ is given by 
\begin{equation*}
K=2mkF+2kR.
\end{equation*}

The cohomology class $K$ is divisible by $2k$. The sphere $S$ sews together with a Seifert surface for the knot $K_2$ to give a surface $C$ in $M$ with $C\cdot R=1$ and $C\cdot F=0$, hence $C\cdot K=2k$. This implies that the divisibility of $K$ is precisely $d=2k$. 

Suppose that $\mathbf{d=2k+1}$ and $\mathbf{n=2m+1}$ are both odd with $k\geq 0$ and $m\geq 1$. We consider the elliptic surface $E(n)$ and do a similar construction. Let $K_1$ be a fibred knot of genus $g_1=2km+k+1$ and do knot surgery along $F$ as above. We get a simply-connected symplectic 4-manifold $M_1$ with canonical class
\begin{align*}
K_{M_1}&=(n-2)F+2g_1F\\
&=(2m+1-2+4km+2k+2)F\\
&=(4km+2k+2m+1)F\\
&=(2m+1)(2k+1)F.
\end{align*}
Next we consider a fibred knot $K_2$ of genus $g_2=2k+1$ and do knot surgery along the rim torus $R$. The result is a simply-connected symplectic 4-manifold $M$ homeomorphic to $E(n)$ with canonical class  
\begin{equation*}
K=(2m+1)(2k+1)F+2(2k+1)R.
\end{equation*}
The class $K$ is divisible by $(2k+1)$. The same argument as above shows that there exists a surface $C$ in $M$ with $C\cdot K=2(2k+1)$. We claim that the divisibility of $K$ is precisely $(2k+1)$: This follows because $M$ is still homeomorphic to $E(n)$ by the knot surgery construction. Since $n$ is odd, the manifold $M$ is not spin and hence $2$ does not divide $K$ (an explicit surface with odd intersection number can be constructed from a section of $E(n)$ and a Seifert surface for the knot $K_1$. This surface has self-intersection number $-n$ and intersection number $(2m+1)(2k+1)$ with $K$.)

To cover the remaining case $m=0$ (corresponding to $n=1$) we can do knot surgery on the elliptic surface $E(1)$ along a general fibre $F$ with a knot $K_1$ of genus $g_1=k+1$. The resulting manifold $M_1$ has canonical class 
\begin{equation*}
K_{M_1}=-F+(2k+2)F=(2k+1)F.
\end{equation*}
   
Suppose that $\mathbf{d=2k+1}$ is odd and $\mathbf{n=2m}$ is even with $k\geq 0$ and $m\geq 1$. We consider the elliptic surface $E(n)$ and first perform a logarithmic transformation along $F$ of index $2$. Let $f$ denote the multiple fibre such that $F$ is homologous to $2f$. There exists a 2-sphere in $E(n)_2$ which intersects $f$ in a single point (for a proof see the following lemma). In particular, the complement of $f$ in $E(n)_2$ is simply-connected. The canonical class of $E(n)_2=E(n)_{2,1}$ is given by
\begin{equation*}
K=(2n-3)f.
\end{equation*}
We can assume that the torus $f$ is symplectic (e.g.~by considering the logarithmic transformation to be done on the complex algebraic surface $E(n)$, resulting in the complex algebraic surface $E(n)_2$). Let $K_1$ be a fibred knot of genus $g_1=4km+k+2$ and do knot surgery along $f$ with $K_1$ as above. The result is a simply-connected symplectic 4-manifold homeomorphic to $E(n)_2$. The canonical class is given by
\begin{align*}
K_{M_1}&=(2n-3)f+2g_1f\\
&=(4m-3+8km+2k+4)f\\
&=(8km+4m+2k+1)f\\
&=(4m+1)(2k+1)f.
\end{align*}
We now consider a fibred knot $K_2$ of genus $g_2=2k+1$ and do knot surgery along the rim torus $R$. We get a simply-connected symplectic 4-manifold $M$ homeomorphic to $E(n)_2$ with canonical class
\begin{equation*}
K=(4m+1)(2k+1)f+2(2k+1)R.
\end{equation*}
A similar argument as above shows that the divisibility of $K$ is $d=2k+1$.
\end{proof}
\begin{lem} Let $p\geq 1$ be an integer and $f$ the multiple fibre in $E(n)_p$. Then there exists a sphere in $E(n)_p$ which intersects $f$ transversely in one point. 
\end{lem}
\begin{proof} We can think of the logarithmic transformation as gluing $T^2\times D^2$ into $E(n)\setminus \mbox{int}\,\nu F$ by a certain diffeomorphism $\phi\colon T^2\times S^1\rightarrow \partial \nu F$. The fibre $f$ corresponds to $T^2\times \{0\}$. Consider a disk of the form $\{*\}\times D^2$. It intersects $f$ once and its boundary maps under $\phi$ to a certain simple closed curve on $\partial \nu F$. Since $E(n)\setminus \mbox{int}\,\nu F$ is simply-connected, this curve bounds a disk in $E(n)\setminus \mbox{int}\,\nu F$. The union of this disk and the disk $\{*\}\times D^2$ is a sphere in $E(n)_p$ which intersects $f$ transversely once. 
\end{proof}
\begin{rem} Under the assumptions of Theorem \ref {elliptic div} it is possible to construct infinitely many homeomorphic but pairwise non-diffeomorphic symplectic homotopy elliptic surfaces $(M_r)_{r\in\mathbb{N}}$ with $\chi_h(M_r)=n$, whose canonical classes all have divisibility equal to $d$. This follows because we can vary in each case the knot $K_1$ and its genus $g_1$ without changing the divisibility of the canonical class. The claim then follows by the formula for the Seiberg-Witten invariants of knot surgery manifolds \cite{FSknot}.
\end{rem}

\section{Generalized knot surgery}\label{sect gen knot surg}

Symplectic manifolds with $c_1^2>0$ and divisible canonical class can be constructed with a version of knot surgery for higher genus surfaces described in \cite{FSfam}. Let $K=K_h$ denote the $(2h+1,-2)$-torus knot, which is a fibred knot of genus $h$. Consider the manifold $M_{K}\times S^1$ from the knot surgery construction, cf.~Section \ref{sect knot surg}. This manifold has the structure of a $\Sigma_h$--bundle over $T^2$:
\begin{equation*}
\begin{CD}
M_{K}\times S^1 @<<< \Sigma_{h}\\
@VVV \\
T^2
\end{CD}
\end{equation*}
We denote a fibre of this bundle by $\Sigma_F$. The fibration defines a trivialization of the normal bundle $\nu\Sigma_F$. We form $g$ consecutive generalized fibre sums along the fibres $\Sigma_F$ to get
\begin{equation*}
Y_{g,h}=(M_K\times S^1)\#_{\Sigma_F=\Sigma_F}\#\dotsc\#_{\Sigma_F=\Sigma_F}(M_K\times S^1).
\end{equation*}
The gluing diffeomorphism is chosen such that it identifies the $\Sigma_h$ fibres in the boundary of the tubular neighbourhoods. This implies that $Y_{g,h}$ is a $\Sigma_h$--bundle over $\Sigma_g$:
\begin{equation*}
\begin{CD}
Y_{g,h} @<<< \Sigma_{h}\\
@VVV \\
\Sigma_g
\end{CD}
\end{equation*}
We denote the fibre again by $\Sigma_F$. The fibre bundle has a section $\Sigma_S$ sewed together from $g$ torus sections of $M_K\times S^1$. Since the knot $K$ is a fibred knot, the manifold $M_K\times S^1$ admits a symplectic structure such that the fibre and the section are symplectic. By the Gompf construction this is then also true for $Y_{g,h}$.   

The invariants of the 4-manifold $Y_{g,h}$ can be calculated by standard formulas, cf.~\cite[Lemma 2.4]{Pa2}:
\begin{align*}
c_1^2(Y_{g,h})&=8(g-1)(h-1)\\
e(Y_{g,h})&=4(g-1)(h-1)\\
\sigma(Y_{g,h})&=0.
\end{align*}
By induction on $g$ one can show that the fundamental group $\pi_1(Y_{g,h})$ is normally generated by the image of $\pi_1(\Sigma_S)$ under inclusion \cite[Proposition 2]{FSfam}. This fact, together with the exact sequence 
\begin{equation*}
H_1(\Sigma_F)\rightarrow H_1(Y_{g,h})\rightarrow H_1(\Sigma_g)\rightarrow 0
\end{equation*}
coming from the long exact homotopy sequence for the fibration $\Sigma_F\rightarrow Y_{g,h}\rightarrow \Sigma_g$ by abelianization, shows that the inclusion $\Sigma_S\rightarrow Y_{g,h}$ induces an isomorphism on $H_1$ and the inclusion $\Sigma_F\rightarrow Y_{g,h}$ induces the zero map. In particular, the homology group $H_1(Y_{g,h};\mathbb{Z})$ is free abelian of rank
\begin{equation*}
b_1(Y_{g,h})=gb_1(M_K\times S^1)=2g.
\end{equation*}
This implies with the formula for the Euler characteristic above
\begin{equation*}
b_2(Y_{g,h})=4h(g-1)+2.
\end{equation*}
The summand $4h(g-1)$ results from $2h$ {\em split classes} (or {\em vanishing classes}) together with $2h$ dual rim tori which are created in each fibre sum. The split classes are formed in the following way: In each fibre sum, the interior of a tubular neighbourhood $\nu\Sigma_F$ of a fibre on each side of the sum is deleted and the boundaries $\partial\nu\Sigma_F$ glued together such that the fibres inside the boundary get identified pairwise. Since the inclusion of the fibre $\Sigma_F$ into $M_K\times S^1$ induces the zero map on first homology, the $2h$ generators of $H_1(\Sigma_h)$, where $\Sigma_h$ is considered as a fibre in $\partial\nu\Sigma_F$, bound surfaces in $M_K\times S^1$ minus the interior of the tubular neighbourhood $\nu\Sigma_F$. The split classes arise from sewing together surfaces bounding corresponding generators on each side of the fibre sum. Fintushel and Stern show that in the case above there exists a basis for the group of split classes consisting of $2h(g-1)$ disjoint surfaces of genus $2$ and self-intersection $2$. This implies
\begin{equation*}
H^2(Y_{g,h};\mathbb{Z})=2h(g-1)\left(\begin{array}{cc} 2 &1 \\ 1& 0 \\ \end{array}\right)\oplus \left(\begin{array}{cc} 0 &1 \\ 1& 0 \\ \end{array}\right),
\end{equation*}
where the last summand is the intersection form on $(\mathbb{Z}\Sigma_S\oplus\mathbb{Z}\Sigma_F)$. They also show that the canonical class of $Y_{g,h}$ is given by
\begin{equation*}
K_Y=(2h-2)\Sigma_S+(2g-2)\Sigma_F,
\end{equation*}
where $\Sigma_S$ and $\Sigma_F$ are oriented by the symplectic form.

Suppose that $M$ is a closed symplectic 4-manifold which contains a symplectic surface $\Sigma_M$ of genus $g$ and self-intersection $0$, oriented by the symplectic form and representing an indivisible homology class. We can then form the symplectic generalized fibre sum
\begin{equation*}
X=M\#_{\Sigma_M=\Sigma_S}Y_{g,h}.
\end{equation*}
If the manifolds $M$ and $M\setminus \Sigma_M$ are simply-connected, then $X$ is again simply-connected because the fundamental group of $Y_{g,h}$ is normally generated by the image of $\pi_1(\Sigma_S)$. Since the inclusion of the surface $\Sigma_S$ in $Y_{g,h}$ induces an isomorphism on first homology, it follows by Theorem \ref{thm formula second cohom fibre sum} and Lemma \ref{lem no rim tori X isom M} that
\begin{equation*}
H^2(X;\mathbb{Z})=P(M)\oplus P(Y_{g,h})\oplus (\mathbb{Z}B_X\oplus\mathbb{Z}\Sigma_X).
\end{equation*}
The surface $B_X$ is sewed together from a surface $B_M$ in $M$ with $B_M\Sigma_M=1$ and the fibre $\Sigma_F$ in the manifold $Y_{g,h}$. Since $\Sigma_F^2=0$, the embedding $H^2(M;\mathbb{Z})\rightarrow H^2(X;\mathbb{Z})$ given by equation \eqref{eq monomorph M into X} preserves the intersection form. Therefore we can write
\begin{equation}\label{decomp of H^2 for Y_g,h}
H^2(X;\mathbb{Z})=H^2(M;\mathbb{Z})\oplus P(Y_{g,h})
\end{equation}
with intersection form
\begin{equation*}
Q_X=Q_M\oplus 2h(g-1)\left(\begin{array}{cc} 2 &1 \\ 1& 0 \\ \end{array}\right).
\end{equation*}
The invariants of $X$ are given by
\begin{align*}
c_1^2(X)&=c_1^2(M)+8h(g-1)\\
e(X)&=e(M)+4h(g-1)\\
\sigma(X)&=\sigma(M).
\end{align*}
The canonical class of $X$ can be calculated by Theorem \ref{thm formula can class gen fibre sum} 
\begin{equation}\label{Fint Stern knot Y}
K_X=K_M+2h\Sigma_M,
\end{equation}
where the isomorphism in equation \eqref{decomp of H^2 for Y_g,h} is understood (this formula follows also from the calculation of Seiberg-Witten invariants in \cite{FSfam}). The formula in equation \eqref{Fint Stern knot Y} is a generalization of the formula in \eqref{Fint Stern knot K}. In particular we get:

\begin{prop}\label{knot surgery high genus div can} Let $M$ be a closed, symplectic 4-manifold which contains a symplectic surface $\Sigma_M$ of genus $g>1$ and self-intersection $0$. Suppose that $\pi_1(M)=\pi_1(M\setminus \Sigma_M)=1$ and that the canonical class of $M$ is divisible by $d$. 
\begin{itemize}
\item If $d$ is odd, there exists for every integer $t\geq 1$ a simply-connected symplectic 4-manifold $X$ with invariants
\begin{align*}
c_1^2(X)&=c_1^2(M)+8td(g-1)\\
e(X)&=e(M)+4td(g-1)\\
\sigma(X)&=\sigma(M)
\end{align*}
and canonical class divisible by $d$.
\item If $d$ is even, there exists for every integer $t\geq 1$ a simply-connected symplectic 4-manifold $X$ with invariants
\begin{align*}
c_1^2(X)&=c_1^2(M)+4td(g-1)\\
e(X)&=e(M)+2td(g-1)\\
\sigma(X)&=\sigma(M)
\end{align*}
and canonical class divisible by $d$.
\end{itemize}
\end{prop}
This follows from the construction above by taking the genus of the torus knot $h=td$ if $d$ is odd and $h=\tfrac{1}{2}td$ if $d$ is even. Hence if a symplectic surface $\Sigma_M$ of genus $g>1$ and self-intersection $0$ exists in $M$, we can raise $c_1^2$ without changing the signature or the divisibility of the canonical class. 

\section{Spin symplectic 4-manifolds with $c_1^2>0$ and negative signature}\label{sect geo spin neg sigma c pos}

We can apply the construction from the previous section to the symplectic homotopy elliptic surfaces constructed in Theorem \ref{elliptic div}. In this section we consider the case of even divisibility $d$ and in the following section the case of odd $d$. 

Recall that in the first case in the proof of Theorem \ref{elliptic div} we constructed a simply-connected symplectic 4-manifold $M$ from the elliptic surface $E(2m)$ by doing knot surgery along a general fibre $F$ with a fibred knot $K_1$ of genus $g_1=(k-1)m+1$ and a further knot surgery along a rim torus $R$ with a fibred knot $K_2$ of genus $g_2=k$. Here $2m\geq 2$ and $d=2k\geq 2$ are arbitrary even integers. The canonical class is given by
\begin{equation*}
K_M=2mkF+2kR=mdF+dR.
\end{equation*}
The manifold $M$ is still homeomorphic to $E(2m)$. There exists an embedded 2-sphere $S$ in $E(2m)$ of self-intersection $-2$ which intersects the rim torus $R$ once. The sphere $S$ is in a natural way Lagrangian \cite{AMP}. We can assume that $S$ is disjoint from the fibre $F$ and by a perturbation of the symplectic structure on $E(2m)$ that the regular fibre $F$, the rim torus $R$ and the dual 2-sphere $S$ are all symplectic and the symplectic form induces a positive volume form on each of them, cf.~the proofs of \cite[Lemma 2.1]{FScan} and \cite[Proposition 3.2]{V2}. 

The 2-sphere $S$ minus a disk sews together with a Seifert surface for $K_2$ to give a symplectic surface $C$ in $M$ of genus $k$ and self-intersection $-2$ which intersects the rim torus $R$ once. By smoothing the double point we get a symplectic surface $\Sigma_M$ in $M$ of genus $g=k+1$ and self-intersection $0$ which represents $C+R$. 

The complement of the surface $\Sigma_M$ in $M$ is simply-connected: This follows because we can assume that $R\cup S$ in the elliptic surface $E(2m)$ is contained in an embedded nucleus $N(2)$, cf.~\cite{GoM, GS} and Example \ref{n-1 Lagrangian triples in E(n)}. Inside the nucleus $N(2)$ there exists a cusp which is homologous to $R$ and disjoint from it. The cusp is still contained in $M$ and intersects the surface $\Sigma_M$ once. Since $M$ is simply-connected and the cusp homeomorphic to $S^2$, the claim $\pi_1(M\setminus \Sigma_M)=1$ follows\footnote{This argument is similar to the argument showing that the complement of a section in $E(n)$ is simply-connected, cf.~\cite[Example 5.2]{Go}.}. 

Let $t\geq 1$ be an arbitrary integer and $K_3$ the $(2h+1,-2)$-torus knot of genus $h=tk$. Consider the generalized fibre sum
\begin{equation*}
X=M\#_{\Sigma_M=\Sigma_S}Y_{g,h}
\end{equation*}
where $g=k+1$. Then $X$ is a simply-connected symplectic 4-manifold with invariants
\begin{align*}
c_1^2(X)&=8tk^2=2td^2\\
e(X)&=24m+4tk^2=24m+td^2\\
\sigma(X)&=-16m.
\end{align*}
The canonical class is given by
\begin{align*}
K_X&=K_M+2tk\Sigma_M\\
&=d(mF+R+t\Sigma_M).
\end{align*}
Hence $K_X$ has divisibility $d$, since the class $mF+R+t\Sigma_M$ has intersection $1$ with $\Sigma_M$. We get:
\begin{thm}\label{thm existence c1 positive d even} Let $d\geq 2$ be an even integer. Then for every pair of positive integers $m,t$ there exists a simply-connected closed spin symplectic 4-manifold $X$ with invariants
\begin{align*}
c_1^2(X)&=2td^2\\
e(X)&=td^2+24m\\
\sigma(X)&=-16m,
\end{align*}
such that the canonical class $K_X$ has divisibility $d$.
\end{thm}
Note that this solves by Lemma \ref{lem 1 geography K} and Rochlin's theorem the existence question for simply-connected 4-manifolds with canonical class divisible by an even integer and negative signature. In particular (for $d=2$), every possible lattice point with $c_1^2>0$ and $\sigma<0$ can be realized by a simply-connected spin symplectic 4-manifold with this construction (the existence of such 4-manifolds has been proved before in \cite{DPZ} in a similar way).

\begin{ex}[Spin homotopy Horikawa surfaces]\label{horikawa with div K} To identify the homeomorphism type of some of the manifolds in Theorem \ref{thm existence c1 positive d even}, let $d=2k$, hence
\begin{align*}
c_1^2(X)&=8tk^2\\
\chi_h(X)&=tk^2+2m.
\end{align*}
We consider the case when the invariants are on the Noether line $c_1^2=2\chi_h-6$: This happens if and only if
\begin{equation*}
6tk^2=4m-6
\end{equation*}
hence $2m=3tk^2+3$, which has a solution if and only if both $t$ and $k$ are odd. Hence for every pair of odd integers $t, k\geq 1$ there exists a simply-connected symplectic 4-manifold $X$ with invariants
\begin{align*}
c_1^2(X)&=8tk^2\\
\chi_h(X)&=4tk^2+3
\end{align*}
such that the divisibility of $K_X$ is $2k$. 

By a construction of Horikawa \cite{HorI} there exists for every odd integer $r\geq 1$ a simply-connected spin complex algebraic surface $M$ on the Noether line with invariants
\begin{align*}
c_1^2(M)&=8r\\
\chi_h(M)&=4r+3.
\end{align*}
See also \cite[Theorem 7.4.20]{GS} where this surface is called $U(3, r+1)$. 

By Freedman's theorem \cite{Freedm} the symplectic 4-manifolds $X$ constructed above for odd parameters $t$ and $k$ are homeomorphic to spin Horikawa surfaces with $r=tk^2$. If $k>1$ and $t$ is arbitrary, the canonical class of $X$ has divisibility $2k>2$. In this case the manifold $X$ cannot be {\em diffeomorphic} to a Horikawa surface: It is known by \cite{HorI} that all Horikawa surfaces $M$ have a fibration in genus $2$ curves, hence by Lemma \ref{lem 2 geography K} the divisibility of $K_M$ is at most $2$ and in the spin case it is equal to $2$. Since Horikawa surfaces are minimal complex surfaces of general type, the claim follows by Proposition \ref{prop no complex surface has two symp with diff div}. 
\end{ex}

\section{Non-spin symplectic 4-manifolds with $c_1^2>0$ and negative signature}\label{sect geo non-spin neg sigma c pos}

In this section we construct some families of simply-connected symplectic 4-manifolds with $c_1^2>0$ such that the divisibility of $K$ is a given odd integer $d>1$. However, we do not have a complete existence result as in Theorem \ref{thm existence c1 positive d even}.

We consider the case that the canonical class $K_X$ is divisible by an odd integer $d$ and the signature $\sigma(X)$ is divisible by $8$. 
\begin{lem} Let $X$ be a closed simply-connected symplectic 4-manifold such that $K_X$ is divisible by an odd integer $d\geq 1$ and $\sigma(X)$ is divisible by $8$. Then $c_1^2(X)$ is divisible by $8d^2$.
\end{lem}
\begin{proof} Suppose that $\sigma(X)=8m$ for some integer $m\in \mathbb{Z}$. Then $b_2^-(X)=b_2^+(X)-8m$ hence $b_2(X)=2b_2^+(X)-8m$. This implies
\begin{equation*}
e(X)=2b_2^+(X)+2-8m.
\end{equation*}
Since $X$ is symplectic, the integer $b_2^+(X)$ is odd, so we can write $b_2^+(X)=2k+1$ for some $k\geq 0$. This implies
\begin{equation*}
e(X)=4k+4-8m,
\end{equation*}
hence $e(X)$ is divisible by $4$. The equation $c_1^2(X)=2e(X)+3\sigma(X)$ shows that $c_1^2(X)$ is divisible by $8$. Since $c_1^2(X)$ is also divisible by the odd integer $d^2$, the claim follows.
\end{proof}

The following theorem covers the case that $K_X$ has odd divisibility and the signature is negative, divisible by $8$ and $\leq -16$:
\begin{thm}\label{odd div sign div 8} Let $d\geq 1$ be an odd integer. Then for every pair of positive integers $n,t$ with $n\geq 2$ there exists a simply-connected closed non-spin symplectic 4-manifold $X$ with invariants
\begin{align*}
c_1^2(X)&=8td^2\\
e(X)&=4td^2+12n\\
\sigma(X)&=-8n
\end{align*}
such that the canonical class $K_X$ has divisibility $d$. 
\end{thm}  
\begin{proof} The proof is similar to the proof of Theorem \ref{thm existence c1 positive d even}. We can write $d=2k+1$ with $k\geq 0$. Suppose that $\mathbf{n=2m+1}$ is odd where $m\geq 1$. In the proof of Theorem \ref{elliptic div} a homotopy elliptic surface $M$ with $\chi_h(M)=n$ was constructed from the elliptic surface $E(n)$ by doing knot surgery along a general fibre $F$ with a fibred knot $K_1$ of genus $g_1=2km+k+1$ and a further knot surgery along a rim torus $R$ with a fibred knot $K_2$ of genus $g_2=2k+1=d$. The canonical class is given by
\begin{align*}
K_M&=(2m+1)(2k+1)F+2(2k+1)R\\
&=(2m+1)dF+2dR.
\end{align*}
There exist a symplectically embedded 2-sphere $S$ in $E(n)$ of self-intersection $-2$ which sews together with a Seifert surface for $K_2$ to give a symplectic surface $C$ in $M$ of genus $d$ and self-intersection $-2$ which intersects the rim torus $R$ once. By smoothing the double point we get a symplectic surface $\Sigma_M$ in $M$ of genus $g=d+1$ and self-intersection $0$ which represents $C+R$. Using a cusp which intersects $\Sigma_M$ once, it follows as above that the complement $M\setminus \Sigma_M$ is simply-connected.

Let $t\geq 1$ be an arbitrary integer and $K_3$ the $(2h+1,-2)$-torus knot of genus $h=td$. We consider the generalized fibre sum
\begin{equation*}
X=M\#_{\Sigma_M=\Sigma_S}Y_{g,h}
\end{equation*}
where $g=d+1$. Then $X$ is a simply-connected symplectic 4-manifold with invariants
\begin{align*}
c_1^2(X)&=8td^2\\
e(X)&=4td^2+12n\\
\sigma(X)&=-8n.
\end{align*}
The canonical class is given by
\begin{align*}
K_X&=K_M+2td\Sigma_M\\
&=d((2m+1)F+2R+2t\Sigma_M).
\end{align*}
Hence $K_X$ has divisibility $d$, since the class $(2m+1)F+2R+2t\Sigma_M$ has intersection $2$ with $\Sigma_M$ and intersection $(2m+1)$ with a surface coming from a section of $E(n)$ and a Seifert surface for $K_1$. 

The case that $\mathbf{n=2m}$ is even where $m\geq 1$ can be proved similarly. By doing a logarithmic transform on the fibre $F$ in $E(n)$ and two further knot surgeries with a fibred knot $K_1$ of genus $g_1=4km+k+2$ on the multiple fibre $f$ and with a fibred knot $K_2$ of genus $g_2=2k+1=d$ along a rim torus $R$, we get a homotopy elliptic surface $M$ with $\chi_h(M)=n$ and canonical class
\begin{equation*}
K_M=(4m+1)df+2dR.
\end{equation*}
The same construction as above yields a simply-connected symplectic 4-manifold $X$ with invariants
\begin{align*}
c_1^2(X)&=8td^2\\
e(X)&=4td^2+12n\\
\sigma(X)&=-8n.
\end{align*}
The canonical class is given by
\begin{align*}
K_X&=K_M+2td\Sigma_M\\
&=d((4m+1)f+2R+2t\Sigma_M).
\end{align*}
Hence $K_X$ has again divisibility $d$. 
\end{proof}

\begin{ex}[Non-spin homotopy Horikawa surfaces] The manifolds in Theorem \ref{odd div sign div 8} have invariants
\begin{align*}
c_1^2(X)&=8td^2\\
\chi_h(X)&=td^2+n.
\end{align*}
Similarly to Example \ref{horikawa with div K} this implies that for every pair of positive integers $d,t\geq 1$ with $d$ odd and $t$ arbitrary there exists a non-spin symplectic homotopy Horikawa surface $X$ on the Noether line $c_1^2=2\chi_h-6$ with invariants
\begin{align*}
c_1^2(X)&=8td^2\\
\chi_h(X)&=4td^2+3,
\end{align*}
whose canonical class has divisibility $d$. Note that for every integer $s\geq 1$ there exists a non-spin complex Horikawa surface $M$ \cite{HorI} with invariants
\begin{align*}
c_1^2(M)&=8s\\
\chi_h(M)&=4s+3.
\end{align*}
If $d>1$ and $t$ an arbitrary integer we get non-spin homotopy Horikawa surfaces with $s=td^2$ whose canonical classes have divisibility $d$. By the same argument as before these 4-manifolds cannot be diffeomorphic to complex Horikawa surfaces.
\end{ex}

With different constructions it is possible to find examples of simply-connected symplectic 4-manifolds with canonical class of odd divisibility, $c_1^2>0$ and signature not divisible by $8$, cf.~\cite[Section VI.2.3]{MHthesis}. However, many cases remain uncovered. For example, we could not answer the following question: 
\begin{quest} For a given odd integer $d>1$ find a simply-connected symplectic 4-manifold $M$ with $c_1^2(M)=d^2$ whose canonical class has divisibility $d$.
\end{quest}
Note that there is a trivial example for $d=3$, namely $\mathbb{C}P^2$.

\section{Construction of inequivalent symplectic structures}\label{sect inequiv symp str}

In this section we prove a result similar to a theorem of Smith \cite[Theorem 1.5]{Smi} which can be used to show that certain 4-manifold $X$ admit inequivalent symplectic structures, where ``equivalence'' is defined in the following way, cf.~\cite{McT}:
\begin{defn} Two symplectic forms on a closed oriented 4-manifold $M$ are called {\em equivalent}, if they can be made identical by a combination of deformations through symplectic forms and orientation preserving self-diffeomorphisms of $M$.
\end{defn}
The canonical classes of equivalent symplectic forms have the same (maximal) divisibility as elements of $H^2(M;\mathbb{Z})$. This follows because deformations do not change the canonical class and the application of an orientation preserving self-diffeomorphism does not change the divisibility.

We will use the following lemma.
\begin{lem}\label{lemma -omega -K} Let $(M,\omega)$ be a symplectic 4-manifold with canonical class $K$. Then the symplectic structure $-\omega$ has canonical class $-K$.
\end{lem}
\begin{proof} Let $J$ be an almost complex structure on $M$, compatible with $\omega$. Then $-J$ is an almost complex structure compatible with $-\omega$. The complex vector bundle $(TM,-J)$ is the conjugate bundle to $(TM,J)$. By \cite{MSt} this implies that $c_1(TM,-J)=-c_1(TM,J)$. Since the canonical class is minus the first Chern class of the tangent bundle, the claim follows. 
\end{proof} 
Let $M_K\times S^1$ be a 4-manifold used in knot surgery where $K$ is a fibred knot of genus $h$. Let $T_K$ be a section of the fibre bundle
\begin{equation*}
\begin{CD}
M_K\times S^1 @<<< \Sigma_{h}\\
@VVV \\
T^2
\end{CD}
\end{equation*}
and $B_K$ a fibre. We fix an orientation on $T_K$ and choose the orientation on $B_K$ such that $T_K\cdot B_K=+1$. There exist symplectic structures on $M_K\times S^1$ such that both the fibre and the section are symplectic. We can choose such a symplectic structure $\omega^+$ which restricts to both $T_K$ and $B_K$ as a positive volume form with respect to the orientations. It has canonical class
\begin{equation*}
K^+=(2h-2)T_K
\end{equation*}
by the adjunction formula. We also define the symplectic form $\omega^-=-\omega^+$. It restricts to a negative volume form on $T_K$ and $B_K$. The canonical class of this symplectic structure is
\begin{equation*}
K^-=-(2h-2)T_K.
\end{equation*} 
Let $X$ be a closed oriented 4-manifold with torsion free cohomology which contains an embedded oriented torus $T_X$ of self-intersection $0$, representing an indivisible homology class. We form the oriented 4-manifold 
\begin{equation*}
X_K=X\#_{T_X=T_K}(M_K\times S^1),
\end{equation*}
by doing the generalized fibre sum along the pair $(T_X,T_K)$ of oriented tori. Suppose that $X$ has a symplectic structure $\omega_X$ such that $T_X$ is symplectic. We consider two cases: If the symplectic form $\omega_X$ restricts to a positive volume form on $T_X$ we can glue this symplectic form to the symplectic form $\omega^+$ on $M_K\times S^1$ to get a symplectic structure $\omega_{X_K}^+$ on $X_K$. The canonical class of this symplectic structure is
\begin{equation*}
K^+_{X_K}=K_X+2hT_X,
\end{equation*}
as seen above, cf.~equation \eqref{Fint Stern knot K}.
\begin{lem}\label{lem knot surgery wrong or} Suppose that $\omega_X$ restricts to a negative volume form on $T_X$. We can glue this symplectic form to the symplectic form $\omega^-$ on $M_K\times S^1$ to get a symplectic structure $\omega_{X_K}^-$ on $X_K$. The canonical class of this symplectic structure is
\begin{equation*}
K^-_{X_K}=K_X-2hT_X.
\end{equation*}
\end{lem}
\begin{proof} We use Lemma \ref{lemma -omega -K} twice: The symplectic form $-\omega_X$ restricts to a positive volume form on $T_X$. We can glue this symplectic form to the symplectic form $\omega^+$ on $M_K\times S^1$ which also restricts to a positive volume form on $T_K$. By the standard formula \eqref{Fint Stern knot K} we get for the canonical class of the resulting symplectic form on $X_K$
\begin{equation*}
K=-K_X+2hT_X.
\end{equation*}
The symplectic form $\omega_{X_K}^-$ we want to consider is {\em minus} the symplectic form we have just constructed. Hence its canonical class is $K^-_{X_K}=K_X-2hT_X$.
\end{proof}

\begin{lem}\label{Lagrang tori in M} Let $(M,\omega)$ be a closed symplectic 4-manifold with canonical class $K_M$. Suppose that $M$ contains pairwise disjoint embedded oriented Lagrangian surfaces $T_1,\dotsc,T_{r+1}$ ($r\geq 1$) with the following properties:
\begin{itemize}
\item The classes of the surfaces $T_1,\dotsc,T_r$ are linearly independent in $H_2(M;\mathbb{R})$.  
\item The surface $T_{r+1}$ is homologous to $a_1T_1+\dotsc +a_rT_r$, where all coefficients $a_1,\dotsc,a_r$ are positive integers. 
\end{itemize} 
Then for every non-empty subset $S\subset \{T_1,\dotsc,T_r\}$ there exists a symplectic form $\omega_S$ on $M$ with the following properties: 
\begin{itemize}
\item All surfaces $T_1,\dotsc,T_{r+1}$ are symplectic. 
\item The symplectic form $\omega_S$ induces on the surfaces in $S$ and the surface $T_{r+1}$ a positive volume form and on the remaining surfaces in $\{T_1,\dotsc,T_r\}\setminus S$ a negative volume form.
\end{itemize}
Moreover, the canonical classes of the symplectic structures $\omega_S$ are all equal to $K_M$. We can also assume that any given closed oriented surface in $M$ disjoint from the surfaces $T_1,\dotsc,T_{r+1}$, which is symplectic with respect to $\omega$, stays symplectic for $\omega_S$ with the same sign of the induced volume form.  
\end{lem} 
\begin{proof} The proof is similar to the proof of \cite[Lemma 1.6]{Go}. We can assume that $S=\{T_{s+1},\dotsc,T_r\}$ with $s+1\leq r$. Let
\begin{equation*}
c=\sum_{i=1}^s a_i,\quad c'=\sum_{i=s+1}^{r-1}a_i.
\end{equation*}
Since the classes of the surfaces $T_1,\dotsc,T_r$ are linearly independent in $H_2(M;\mathbb{R})$ and $H^2_{DR}(M)$ is the dual space of $H_2(M;\mathbb{R})$ there exists a closed 2-form $\eta$ on $M$ with the following properties:
\begin{align*}
\int_{T_1}\eta&=-1,\dotsc, \int_{T_s}\eta =-1\\  
\int_{T_{s+1}}\eta&=+1,\dotsc,\int_{T_{r-1}}\eta=+1\\
\int_{T_r}\eta&=c+1\\
\int_{T_{r+1}}\eta&=c'+1.
\end{align*}
Note that we can choose the value of $\eta$ on $T_1,\dotsc,T_r$ arbitrarily. The value on $T_{r+1}$ is then determined by $T_{r+1}=a_1T_1+\dotsc+a_rT_r$. We can choose symplectic forms $\omega_i$ on each $T_i$ such that
\begin{equation*}
\int_{T_i}\omega_i=\int_{T_i}\eta,\quad\mbox{for all $i=1,\dotsc,r+1$}.
\end{equation*}
The symplectic form $\omega_i$ induces on $T_i$ a negative volume form if $i\leq s$ and a positive volume form if $i\geq s+1$. The difference $\omega_i-j_i^*\eta$, where $j_i\colon T_i\rightarrow M$ is the embedding, has vanishing integral and hence is an exact 2-form on $T_i$ of the form $d\alpha_i$. We can extend each $\alpha_i$ to a small tubular neighbourhood of $T_i$ in $M$, cut it off differentiably in a slightly larger tubular neighbourhood and extend by $0$ to all of $M$. We can do this such that the tubular neighbourhoods of $T_1,\dotsc,T_{r+1}$ are pairwise disjoint. Define the closed 2-form $\eta'=\eta+\sum_{i=1}^{r+1}d\alpha_i$ on $M$. Then
\begin{equation*}
j_i^*\eta'=j_i^*\eta+d\alpha_i=\omega_i.
\end{equation*}
The closed 2-form $\omega'=\omega+t\eta'$ is for small values of $t$ symplectic. Since the surfaces $T_i$ are Lagrangian, we have $j_i^*\omega=0$ and hence  $j_i^*\omega'=t\omega_i$. This implies that $\omega'$ is for small values $t>0$ a symplectic form on $M$ which induces a volume form on $T_i$ of the same sign as $\omega_i$ for all $i=1,\dotsc,r+1$. The claim about the canonical class follows because the symplectic structures $\omega_S$ are constructed by a deformation of $\omega$. We can also choose $t>0$ small enough such that $\omega'$ still restricts to a symplectic form on any given symplectic surface disjoint from the tori without changing the sign of the induced volume form on this surface.
\end{proof}

This construction will be used in the following way: Suppose that $(V_1,\omega_1)$ and $(V_2,\omega_2)$ are symplectic 4-manifolds, such that $V_1$ contains an embedded Lagrangian torus $T_1$ and $V_2$ contains an embedded symplectic torus $T_2$, both oriented and of self-intersection $0$. Let $W$ denote the smooth oriented 4-manifold $V_1\#_{T_1=T_2}V_2$, obtained as a generalized fibre sum. By the previous lemma, there exist small perturbations of $\omega_1$ to new symplectic forms $\omega_1^+$ and $\omega_1^-$ on the manifold $V_1$, such that the torus $T_1$ becomes symplectic with positive and negative induced volume form, respectively. By the Gompf construction it is then possible to define two symplectic forms on the same oriented 4-manifold $W$:
\begin{enumerate}
\item  The symplectic forms $\omega_1^+$ and $\omega_2$ determine a symplectic form on $W$. 
\item The symplectic forms $\omega_1^-$ and $-\omega_2$ determine a symplectic form on $W$.
\end{enumerate}
Hence the symplectic forms on the first manifold differ only by a small perturbation, while on the second manifold they differ by the sign. Similarly, the canonical classes of both perturbed symplectic forms on $V_1$ are the same, while they differ by the sign on $V_2$. If additional tori exist and suitable fibre sums performed, it is possible to end up with two or more inequivalent symplectic forms on the same 4-manifold, distinguished by the divisibilities of their canonical classes.

To define the configuration of tori we want to consider, recall that the nucleus $N(n)$ is the smooth manifold with boundary defined as a regular neighbourhood of a cusp fibre and a section in the simply-connected elliptic surface $E(n)$, cf.~\cite{Gnuc}. It contains an embedded torus given by a regular fibre homologous to the cusp. It also contains two embedded disks of self-intersection $-1$ which bound vanishing cycles on the torus. The vanishing cycles are the simple-closed loops given by the factors in $T^2=S^1\times S^1$.

\begin{defn}[Lagrangian triple] Let $(M,\omega)$ be a symplectic 4-manifold. Given an integer $a\geq 1$, a {\em Lagrangian triple} consists of three pairwise disjoint oriented Lagrangian tori $T_1,T_2,R$ embedded in $M$ with the following properties:
\begin{itemize}
\item All three tori have self-intersection $0$ and represent indivisible classes in integral homology.
\item $T_1,T_2$ are linearly independent over $\mathbb{Q}$ and $R$ is homologous to $aT_1+T_2$.
\item There exists an embedded nucleus $N(2)\subset M$ which contains the torus $R$, corresponding to a general fibre. Let $S$ denote the 2-sphere in $N(2)$ of self-intersection $-2$, corresponding to a section. In addition to intersecting $R$, this sphere is assumed to intersect $T_2$ transversely once.
\item The torus $T_1$ is disjoint from the nucleus $N(2)$ above and there exists an embedded 2-sphere $S_1$ in $M$, also disjoint from $N(2)$, which intersects $T_1$ transversely and positively once.
\end{itemize}
\end{defn}
\begin{figure}[htbp]
	\centering
		\includegraphics[width=12cm]{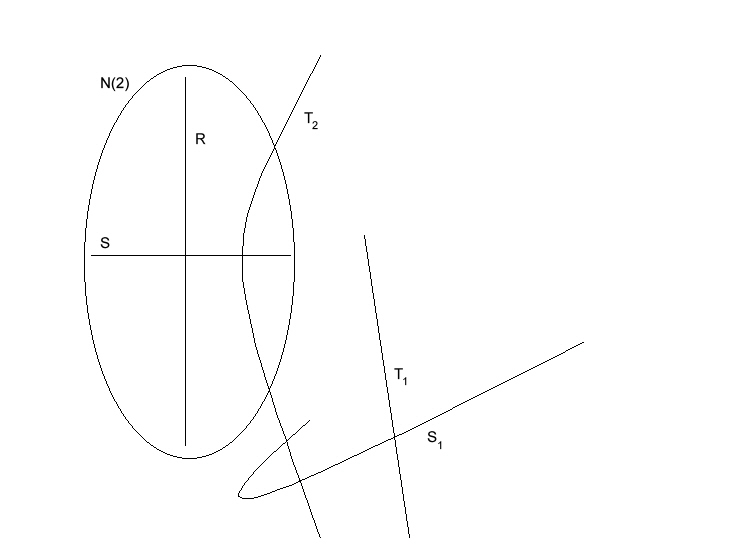}
		\caption{Lagrangian triple}
		\label{figure Lag triple}
\end{figure}
See Figure \ref{figure Lag triple}. The assumptions imply that 
\begin{equation*}
S_1T_2=S_1(R-aT_1)=-a.
\end{equation*}

\begin{ex}\label{n-1 Lagrangian triples in E(n)} Let $M$ be the elliptic surface $E(n)$ with $n\geq 2$. In this example we show that $E(n)$ contains $n-1$ disjoint Lagrangian triples $(T_1^i, T_2^i, R^i)$ as above, where $R^i$ is homologous to $a_iT_1^i+T_2^i$, for $i=1,\dotsc,n-1$. The integers $a_i>0$ can be chosen arbitrarily and for each triple independently. In this case both $T_1^i$ and $R^i$ are contained in disjoint embedded nuclei $N(2)$. Together with their dual 2-spheres they realize $2(n-1)$ $H$-summands in the intersection form of $E(n)$. In particular, the tori in different triples are linearly independent. We can also achieve that all Lagrangian tori and the 2-spheres that intersect them once are disjoint from the nucleus $N(n)\subset E(n)$, defined as a regular neighbourhood of a cusp fibre and a section in $E(n)$.

The construction is similar to \cite[Section 2]{GoM} and is done by induction. Suppose the Lagrangian triples are already constructed for $E(n)$ and consider a splitting of $E(n+1)$ as a fibre sum $E(n+1)=E(n)\#_{F=F}E(1)$ along general fibres $F$. We choose fibred tubular neighbourhoods for the general fibres in $E(n)$ and $E(1)$. The boundary of $E(1)\setminus \mbox{int}\,\nu F$ is diffeomorphic to $F\times S^1$. Let $\gamma_1,\gamma_2$ be two simple closed loops spanning the torus $F$ and $m$ the meridian to $F$, spanning the remaining $S^1$ factor. We consider the following three tori
\begin{align*}
V_0&=\gamma_1\times \gamma_2\\
V_1&=\gamma_1\times m\\
V_2&=\gamma_2\times m.
\end{align*}
The tori are made disjoint by pushing them inside a collar of the boundary into the interior of $E(1)\setminus \mbox{int}\,\nu F$ such that $V_2$ is the innermost and $V_0$ the outermost (closest to the boundary). The torus $V_0$ can be assumed symplectic, while $V_1$ and $V_2$ are rim tori which can be assumed Lagrangian. Similarly the boundary of $E(n)\setminus \mbox{int}\,\nu F$ is diffeomorphic to $F\times S^1$, where $F$ is spanned by the circles $\gamma_1,\gamma_2$ and $S^1$ by the circle $m$ and corresponding circles get identified in the gluing of the fibre sum. In the interior of $E(n)\setminus \mbox{int}\,\nu F$ we consider three tori $V_0,V_1,V_2$ as above which get identified with the corresponding tori on the $E(1)$ side in the gluing. On the $E(n)$ side, the torus $V_0$ is the innermost and $V_2$ the outermost.

We can choose elliptic fibrations such that near the general fibre $F$ there exist two cusp fibres in $E(1)$ and three cusp fibres in $E(n)$. This is possible because $E(m)$ has an elliptic fibration with $6m$ cusp fibres for all $m$, cf.~\cite[Corollary 7.3.23]{GS}. The corresponding vanishing disks can be assumed pairwise disjoint. We can also choose three disjoint sections for the elliptic fibration on $E(1)$ and one section for $E(n)$.

The nuclei can now be defined as follows: The nucleus $N(n+1)$ containing $V_0$ has a dual $-(n+1)$-sphere sewed together from sections on each side of the fibre sum. The vanishing disks for $V_0$ come from the first cusp in $E(n)$. The nucleus $N(2)$ containing $V_1$ has a dual $-2$-sphere sewed together from two vanishing cycles parallel to $\gamma_2$ coming from the first cusp in $E(1)$ and the second cusp in $E(n)$. The vanishing disks for $V_1$ come from the second section of $E(1)$ and from the vanishing cycle parallel to $\gamma_1$ of the second cusp in $E(n)$. The nucleus $N(2)$ containing $V_2$ has a dual $-2$-sphere sewed together from two vanishing cycles parallel to $\gamma_1$ coming from the second cusp in $E(1)$ and the third cusp in $E(n)$. The vanishing disks for $V_1$ come from the third section of $E(1)$ and from the vanishing cycle parallel to $\gamma_2$ of the second cusp in $E(1)$.  

To define the Lagrangian triple $(T_1, T_2, R)$ let $T_1=V_1$ and $R=V_2$. Denote by $c_a\colon S^1\rightarrow F=S^1\times S^1$ the embedded curve given by the $(-a,1)$-torus knot and let $T_2$ denote the Lagrangian rim torus 
\begin{equation*}
T_2=c_a\times m
\end{equation*}
in the collar above. Then $T_2$ represents the class $-aT_1+R$, hence $R=aT_1+T_2$. The torus $T_2$ has one positive transverse intersection with the sphere in the nucleus containing $R$ and $a$ negative transverse intersections with the sphere in the nucleus containing $T_1$. 
\end{ex}

\begin{rem} To find more general examples of symplectic 4-manifolds containing Lagrangian triples, suppose that $Y$ is an arbitrary closed symplectic 4-manifold which contains an embedded symplectic torus $T_Y$ of self-intersection $0$, representing an indivisible class. Then the symplectic generalized fibre sum $Y\#_{T_Y=F}E(n)$ also contains $n-1$ Lagrangian triples. 
\end{rem}

Suppose $(M,\omega)$ is a simply-connected symplectic 4-manifold which contains a Lagrangian triple $T_1, T_2, R$. Let $K_1$ and $K_2$ be fibred knots of genera $h_1, h_2$ to be chosen later. Consider the associated oriented 4-manifolds $M_{K_i}\times S^1$ as in the knot surgery construction and denote sections of the fibre bundles
\begin{equation*}
\begin{CD}
M_{K_i}\times S^1 @<<< \Sigma_{h_i}\\
@VVV \\
T^2
\end{CD}
\end{equation*}
by $T_{K_i}$, which are tori of self-intersection $0$. We choose an orientation on each torus $T_{K_i}$. Note that the Lagrangian tori $T_1,T_2$ in $M$ are oriented {\it a priori}. 

We construct a smooth oriented 4-manifold $X$ in three steps as follows: For an integer $m\geq 1$ consider the elliptic surface $E(m)$ and denote an oriented general fibre by $F$. Let $M_0$ denote the smooth generalized fibre sum
\begin{equation*}
M_0=E(m)\#_{F=R}M.
\end{equation*}
The gluing diffeomorphism is chosen in the following way: the push-offs $R'$ and $F'$ into the boundary of the tubular neighbourhoods $\nu R$ and $\nu F$ each contain a pair of vanishing cycles. We choose the gluing such that the push-offs and the vanishing cycles get identified. The corresponding vanishing disks then sew together pairwise to give two embedded spheres of self-intersection $-2$ in $M_0$, which can be assumed disjoint by choosing two different push-offs given by the same trivializations.

Denote the torus in $M_0$ coming from the push-off $R'$ by $R_0$. There exist two disjoint tori in $M_0$, still denoted by $T_1, T_2$, such that $R_0$ is homologous to $aT_1+T_2$ in $M_0$. In the second step of the construction we do a knot surgery with the fibred knot $K_1$ along the torus $T_1$ in $M_0$ to get the oriented 4-manifold
\begin{equation*}
M_1=M_0\#_{T_1=T_{K_1}}(M_{K_1}\times S^1).
\end{equation*}
The manifold $M_1$ contains a torus which we still denote by $T_2$. We do a knot surgery with the fibred knot $K_2$ along the torus $T_2$ to get the oriented 4-manifold
\begin{equation*}
X=M_1\#_{T_2=T_{K_2}}(M_{K_2}\times S^1). 
\end{equation*}

\begin{lem} The closed oriented 4-manifold 
\begin{equation*}
X=E(m)\#_{F=R}M\#_{T_1=T_{K_1}}(M_{K_1}\times S^1)\#_{T_2=T_{K_2}}(M_{K_2}\times S^1)
\end{equation*}
is simply-connected.
\end{lem}
\begin{proof} The existence of the sphere $S$ shows that $M\setminus R$ is simply-connected. Since $E(m)\setminus F$ is simply-connected, it follows that $M_0$ is simply-connected. 

The sphere $S$ and a section for the elliptic fibration on $E(m)$ sew together to give an embedded sphere $S_2$ in $M_0$ of self-intersection $-(m+2)$. The sphere $S_1$ in $M$ is disjoint from $R$ and hence is still contained in $M_0$. These spheres have the following intersections:
\begin{itemize}
\item The sphere $S_1$ intersects $T_1$ transversely once, has intersection $-a$ with $T_2$ and is disjoint from $R_0$.
\item The sphere $S_2$ intersects $R_0$ and $T_2$ transversely once and is disjoint from $T_1$. 
\end{itemize}
The sphere $S_1$ shows that $M_0\setminus T_1$ is simply-connected, hence $M_1$ is simply-connected. The sphere $S_2$ in $M_0$ is disjoint from $T_1$, hence it is still contained in $M_1$ and intersects $T_2$ once. By the same argument this shows that the manifold $X$ is simply-connected. 
\end{proof}

We define two symplectic forms $\omega_X^+$ and $\omega_X^-$ on $X$: By Lemma \ref{Lagrang tori in M} there exist two symplectic structures $\omega_+, \omega_-$ on $M$ with the same canonical class $K_M$ as $\omega$ such that
\begin{itemize}
\item The tori $T_1,T_2,R$ are symplectic with respect to both symplectic forms. 
\item The form $\omega_+$ induces on $T_1,T_2$ and $R$ a positive volume form.
\item The form $\omega_-$ induces on $T_1$ a negative volume form and on $T_2$ and $R$ a positive volume form.
\end{itemize}
We can also achieve that the sphere $S$ is symplectic with positive volume form in both cases. 

On the elliptic surface $E(m)$ we can choose a symplectic (K\"ahler) form $\omega_E$ which restricts to a positive volume form on the oriented fibre $F$. It has canonical class
\begin{equation*}
K_{E}=(m-2)F.
\end{equation*}
We can glue both symplectic forms $\omega_+$ and $\omega_-$ on $M$ to the symplectic form $\omega_E$ on $E(m)$ to get symplectic forms $\omega_0^+$ and $\omega_0^-$ on the 4-manifold $M_0$. The canonical class for both symplectic forms on $M_0$ is given by
\begin{equation*}
K_{M_0}=K_M+mR_0,
\end{equation*}
cf.~\cite[Proof of Lemma 2.2]{FScan}. Since rim tori exist in this fibre sum, Theorem Theorem \ref{thm formula can class gen fibre sum} cannot be applied directly. However, the formula remains correct, because rim tori do not contribute in this case, cf.~\cite[Section V.6.1]{MHthesis} for details.

We want to extend the symplectic forms to the 4-manifold $X$: First we choose in each fibre bundle $M_{K_i}\times S^1$ a fibre $B_{K_i}$ and orient the surface $B_{K_i}$ such that $T_{K_i}\cdot B_{K_i}=+1$ with the chosen orientation on $T_{K_i}$. There exist symplectic structures on the closed 4-manifolds $M_{K_i}\times S^1$ such that both the section and the fibre are symplectic. On $M_{K_1}\times S^1$ we choose two symplectic forms $\omega_1^\pm$: The form $\omega_1^+$ induces a positive volume form on both $T_{K_1}$ and $B_{K_1}$. It has canonical class
\begin{equation*}
K_1^+=(2h_1-2)T_{K_1}.
\end{equation*}
The form $\omega_1^-$ is given by $-\omega_1^+$. It induces a negative volume form on both $T_{K_1}$ and $B_{K_1}$ and has canonical class
\begin{equation*}
K_1^-=-(2h_1-2)T_{K_1}.
\end{equation*}

On the manifold $M_{K_2}\times S^1$ we only choose a symplectic form $\omega_2$ which induces a positive volume form on $T_{K_2}$ and $B_{K_2}$. The canonical class is given by
\begin{equation*}
K_2=(2h_2-2)T_{K_2}.
\end{equation*}

The oriented torus $T_1$ in $M_0$ is symplectic for both forms $\omega_0^\pm$ constructed above such that $\omega_0^+$ induces a positive volume form and $\omega_0^-$ a negative volume form. Gluing $\omega_0^+$ to $\omega_1^+$ and $\omega_0^-$ to $\omega_1^-$ it follows that the closed oriented 4-manifold $M_1$ has two symplectic structures with canonical classes
\begin{align*}
K_{M_1}^+&= K_M+mR_0+2h_1T_1\\
K_{M_1}^-&= K_M+mR_0-2h_1T_1
\end{align*}
by Lemma \ref{lem knot surgery wrong or}.

The torus $T_2$ can be considered as a symplectic torus in $M_1$ such that both symplectic structures induce positive volume forms, since we can assume that the symplectic forms on $M_1$ are still of the form $\omega_0^+$ and $\omega_0^-$ in a neighbourhood of $T_2$. Hence on the generalized fibre sum $X=M_1\#_{T_2=T_{K_2}}M_{K_2}\times S^1$ we can glue each of the two symplectic forms on $M_1$ to the symplectic form $\omega_2$ on $M_{K_2}\times S^1$. We get two symplectic structures on $X$ with canonical classes
\begin{align*}
K_X^+&= K_M+mR_0+2h_1T_1+2h_2T_2\\
K_X^-&= K_M+mR_0-2h_1T_1+2h_2T_2.
\end{align*}
This can be written using $R_0=aT_1+T_2$ as
\begin{align*}
K_X^+&= K_M+(2h_1+am)T_1+(2h_2+m)T_2\\
K_X^-&= K_M+(-2h_1+am)T_1+(2h_2+m)T_2.
\end{align*}
\begin{thm}\label{triple lag tori construct} Suppose that $(M,\omega)$ is a simply-connected symplectic 4-manifold which contains a Lagrangian triple $T_1,T_2, R$ such that $R$ is homologous to $aT_1+T_2$. Let $m$ be a positive integer and $K_1,K_2$ fibred knots of genus $h_1$ and $h_2$. Then the closed oriented 4-manifold 
\begin{equation*}
X=E(m)\#_{F=R}M\#_{T_1=T_{K_1}}(M_{K_1}\times S^1)\#_{T_2=T_{K_2}}(M_{K_2}\times S^1)
\end{equation*}
is simply-connected and admits two symplectic structures $\omega_X^+,\omega_X^-$ with canonical classes
\begin{align*}
K_X^+&= K_M+(2h_1+am)T_1+(2h_2+m)T_2\\
K_X^-&= K_M+(-2h_1+am)T_1+(2h_2+m)T_2. 
\end{align*} 
\end{thm}
\begin{rem} Instead of doing the generalized fibre sum with $E(m)$ in the first step of the construction, we could also do a knot surgery with a fibred knot $K_0$ of genus $h_0\geq 1$. This has the advantage that both $c_1^2$ and the signature do not change under the construction. However, the sphere $S_2$ in $M_0$ is then replaced by a surface of genus $h_0$, sewed together from the sphere $S$ in $M$ and a Seifert surface for $K_0$. Hence it is no longer clear that $M_1\setminus T_2$ and $X$ are simply-connected.
\end{rem}
The following two surfaces are useful in applications to determine the divisibility of the canonical classes in Theorem \ref{triple lag tori construct}.
\begin{lem}\label{defn surface C2} There exists an oriented surface $C_2$ in $X$ which has intersection $C_2T_2=1$ and is disjoint from $T_1$.
\end{lem}
The surface $C_2$ is sewed together from the sphere $S_2$ and a Seifert surface for $K_2$. 

\begin{lem}\label{lem surface C1} There exists an oriented surface $C_1$ in $X$ which has intersection $C_1T_1=1$ and is disjoint from $T_2$. 
\end{lem}
\begin{proof} The surface $C_1$ can be constructed explicitly as follows: In the nucleus $N(2)\subset M$ containing $R$ we can find a surface of some genus homologous to $aS$ and intersecting both $R$ and $T_2$ in $a$ positive transverse intersections. Tubing this surface to the sphere $S_1$ we get a surface $A$ in $M$ which has intersection number $AT_2=0$ and intersects $T_1$ transversely once. By increasing the genus we can achieve that $A$ is disjoint from $T_2$. The surface $A$ still intersects the torus $R$ in $a$ points. Sewing the surface $A$ to a surface in $E(m)$ homologous to $a$ times a section we get a surface $B$ in $M_0$, disjoint from $T_2$ and intersecting $T_1$ once. Sewing this surface to a Seifert surface for $K_1$ we get a surface $C_1$ in $X$ with $C_1T_1=1$, disjoint from $T_2$. 
\end{proof}

\section{Examples of inequivalent symplectic structures}\label{sect ex inequiv symp str}

This section contains some applications of the construction in Section \ref{sect inequiv symp str}. We start with the following definition:
\begin{defn}[Definition of the set $Q$]\label{defn set Q of divisors} Let $N\geq 0$, $d\geq 1$ be integers and $d_0,\dotsc,d_N$ positive integers dividing $d$, where $d=d_0$. If $d$ is even, assume that all $d_1,\dotsc,d_N$ are even. We define a set $Q$ of positive integers as follows:
\begin{itemize}
\item If $d$ is either odd or not divisible by $4$, let $Q$ be the set consisting of the greatest common divisors of all (non-empty) subsets of $\{d_0,\dotsc,d_N\}$. 
\item If $d$ is divisible by $4$ we can assume by reordering that $d_1,\dotsc,d_s$ are those elements such that $d_i$ is divisible by $4$ while $d_{s+1},\dotsc,d_N$ are those elements such that $d_i$ is not divisible by $4$, where $s\geq 0$ is some integer. Then $Q$ is defined as the set of integers consisting of the greatest common divisors of all (non-empty) subsets of $\{d_0,\dotsc,d_s, 2d_{s+1},\dotsc,2d_N\}$.
\end{itemize}
\end{defn}

We can now formulate the main theorem on the existence of inequivalent symplectic structures on homotopy elliptic surfaces:
\begin{thm}\label{elliptic real q} Let $N,d\geq 1$ be integers and $d_0,\dotsc,d_N$ positive integers dividing $d$ as in Definition \ref{defn set Q of divisors}. Let $Q$ be the associated set of greatest common divisors. Choose an integer $n\geq 3$ as follows:
\begin{itemize}
\item If $d$ is odd, let $n$ be an arbitrary integer with $n\geq 2N+1$. 
\item If $d$ is even, let $n$ be an even integer with $n\geq 3N+1$. 
\end{itemize}
Then there exists a homotopy elliptic surface $W$ with $\chi_h(W)=n$ and the following property: For each integer $q\in Q$ the manifold $W$ admits a symplectic structure whose canonical class $K$ has divisibility equal to $q$. Hence $W$ admits at least $|Q|$ many inequivalent symplectic structures. 
\end{thm}  
\begin{proof} The proof splits into three cases depending on the parity of $d$. In each case we follow the construction in Section \ref{sect inequiv symp str}, starting from the manifold $M=E(l)$, where $l$ is an integer $\geq N+1$. By Example \ref{n-1 Lagrangian triples in E(n)} the manifold $E(l)$ contains $N$ pairwise disjoint Lagrangian triples $T_1^i, T_2^i, R^i$, where $R^i$ is homologous to $a_iT_1^i+T_2^i$, for indices $i=1,\dotsc,N$. The construction is done on each triple separately\footnote{This can be done by a small generalization of Lemma \ref{Lagrang tori in M}, because the construction in the proof of this lemma changes the symplectic structure only in a small neighbourhood of the Lagrangian surfaces.} and involves knot surgeries along $T_1^i$ and $T_2^i$ with fibred knots of genus $h_i$ and $h$, respectively, as well as fibre summing with elliptic surfaces $E(m)$ along the tori $R^i$. The numbers $a_i,h_i, h$ and $m$ will be fixed in each case.

Suppose that {\bf $\mathbf{d}$ is odd}. Then all divisors $d_1,\dotsc,d_N$ are odd. Consider the integers defined by
\begin{align*}
m&=1\\
h&=\tfrac{1}{2}(d-1)\\
a_i&=d+d_i\\
h_i&=\tfrac{1}{2}(d-d_i),\quad\mbox{for $1\leq i\leq N$}.\\
\end{align*}
Let $l$ be an integer $\geq N+1$ and do the construction above, starting from the elliptic surface $E(l)$. We get a (simply-connected) homotopy elliptic surface $X$ with $\chi_h(X)=l+N$. By Theorem \ref{triple lag tori construct} the 4-manifold $X$ has $2^N$ symplectic structures with canonical classes
\begin{equation*}\begin{split}
K_X&=(l-2)F+\sum_{i=1}^{N}\left((\pm 2h_i+a_i)T_1^i+(2h+1)T_2^i\right)\\
&=(l-2)F+\sum_{i=1}^{N}\left((\pm (d-d_i)+d+d_i)T_1^i+dT_2^i\right).
\end{split}\end{equation*} 
Here $F$ denotes the torus in $X$ coming from a general fibre in $E(l)$ and the $\pm$-signs in each summand can be varied independently. We can assume that $F$ is symplectic with positive induced volume form for all $2^N$ symplectic structures on $X$. Consider the even integer $l(d-1)+2$ and let $K$ be a fibred knot of genus $g=\tfrac{1}{2}(l(d-1)+2)$. We do knot surgery with $K$ along the symplectic torus $F$ to get a homotopy elliptic surface $W$ with $\chi_h(W)=l+N$ having $2^N$ symplectic structures whose canonical classes are given by
\begin{equation*}\begin{split}
K_W&=(l-2+2g)F+\sum_{i=1}^{N}\left((\pm (d-d_i)+d+d_i)T_1^i+dT_2^i\right)\\
&=dlF+\sum_{i=1}^{N}\left((\pm (d-d_i)+d+d_i)T_1^i+dT_2^i\right). 
\end{split}\end{equation*}
Suppose that $q\in Q$ is the greatest common divisor of certain elements $\{d_i\}_{i\in I}$, where $I$ is a non-empty subset of $\{0,\dotsc,N\}$. Let $J$ be the complement of $I$ in $\{0,\dotsc,N\}$. Choosing the minus sign for each $i$ in $I$ and the plus sign for each $i$ in $J$ defines a symplectic structure $\omega_I$ on $W$ with canonical class given by
\begin{equation*}
K_W = dlF+\sum_{i\in I}(2d_iT_1^i+dT_2^i)+\sum_{j\in J}(2dT_1^j+dT_2^j).
\end{equation*}
We claim that the divisibility of $K_W$ is equal to $q$. Since $q$ divides $d$ and all integers $d_i$ for $i\in I$, the class $K_W$ is divisible by $q$. Considering the surfaces from Lemma \ref{defn surface C2} and Lemma \ref{lem surface C1} for each Lagrangian triple separately implies that every number that divides $K_W$ is odd (since it divides $d$) and a common divisor of all $d_i$ with indices $i\in I$. This proves the claim.

Suppose that {\bf $\mathbf{d}$ is even but not divisible by $\mathbf{4}$}. We can write $d=2k$ and $d_i=2k_i$ for all $i=1,\dotsc,N$. The assumption implies that all integers $k, k_i$ are odd. Consider the integers defined by
\begin{align*}
m&=2\\
h&=k-1\\
a_i&=\tfrac{1}{2}(k+k_i)\\
h_i&=\tfrac{1}{2}(k-k_i).
\end{align*}
Let $l$ be an even integer $\geq N+1$ and consider the construction above, starting from $E(l)$. We get a homotopy elliptic surface $X$ with $\chi_h(X)=l+2N$. The 4-manifold $X$ has $2^N$ symplectic structures with canonical classes
\begin{equation*}\begin{split}
K_X&=(l-2)F+\sum_{i=1}^{N}\left((\pm 2h_i+2a_i)T_1^i+(2h+2)T_2^i\right)\\
&=(l-2)F+\sum_{i=1}^{N}\left((\pm (k-k_i)+k+k_i)T_1^i+dT_2^i\right).
\end{split}\end{equation*}   
Consider a fibred knot $K$ of genus $g=\tfrac{1}{2}(l(d-1)+2)$ (note that $l$ is even). Doing knot surgery with $K$ along the symplectic torus $F$ in $X$ we get a homotopy elliptic surface $W$ with $\chi_h(W)=l+2N$ having $2^N$ symplectic structures whose canonical classes are
\begin{equation*}\begin{split}
K_W&=(l-2+2g)F+\sum_{i=1}^{N}\left((\pm (k-k_i)+k+k_i)T_1^i+dT_2^i\right)\\
&=dlF+\sum_{i=1}^{N}\left((\pm (k-k_i)+k+k_i)T_1^i+dT_2^i\right). 
\end{split}\end{equation*}
Let $q\in Q$ be the greatest common divisor of elements $d_i$ where $i\in I$ for some non-empty index set $I$ with complement $J$ in $\{0,\dotsc,N\}$. Choosing the plus and minus signs as before, we get a symplectic structure $\omega_I$ on $W$ with canonical class 
\begin{equation}\label{canonical class of W pos c12 d not div 4}
K_W=dlF+\sum_{i\in I}(d_iT_1^i+dT_2^i)+\sum_{j\in J}(dT_1^i+dT_2^i).
\end{equation}
As above, it follows that the canonical class of $\omega_I$ has divisibility equal to $q$. 

Finally we consider the case that {\bf $\mathbf{d}$ is divisible by $4$}. We can write $d=2k$ and $d_i=2k_i$ for all $i=1,\dotsc,N$. We can assume that the divisors are ordered as in Definition \ref{defn set Q of divisors}, i.e.~$d_1,\dotsc,d_s$ are those elements such that $d_i$ is divisible by $4$ while $d_{s+1},\dotsc,d_N$ are those elements such that $d_i$ is not divisible by $4$. This is equivalent to $k_1,\dotsc,k_s$ being even and $k_{s+1},\dotsc,k_N$ odd. Consider the integers defined by
\begin{align*}
a_i&=\tfrac{1}{2}(k+k_i)\\
h_i&=\tfrac{1}{2}(k-k_i),\\
\end{align*}
for $i=1,\dotsc,s$, and
\begin{align*}
a_i&=\tfrac{1}{2}(k+2k_i)\\
h_i&=\tfrac{1}{2}(k-2k_i)\\
\end{align*}
for $i=s+1,\dotsc,N$. We also define 
\begin{align*}
m&=2\\
h&=k-1.
\end{align*}
Let $l$ be an even integer $\geq N+1$. We consider the same construction as above starting from $E(l)$ to get a homotopy elliptic surface $X$ with $\chi_h(X)=l+2N$ that has $2^{N}$ symplectic structures with canonical classes given by the formula
\begin{equation*}
\begin{split}
K_X&=(l-2)F+\sum_{i=1}^{N}\left((\pm 2h_i+2a_i)T_1^i+(2h+2)T_2^i\right)\\
&=(l-2)F+\sum_{i=1}^{s}\left((\pm (k-k_i)+k+k_i)T_1^i+dT_2^i\right)\\
&\quad +\sum_{i=s+1}^{N}\left((\pm (k-2k_i)+k+2k_i)T_1^i+dT_2^i\right).
\end{split}
\end{equation*}     
We then do knot surgery with a fibred knot $K$ of genus $g=\tfrac{1}{2}(l(d-1)+2)$ along the symplectic torus $F$ in $X$ to get a homotopy elliptic surface $W$ with $\chi_h(W)=l+2N$ having $2^N$ symplectic structures whose canonical classes are
\begin{equation}\label{canonical class of W pos c12 d does div 4}
\begin{split}
K_W&=(l-2+2g)F+\sum_{i=1}^{N}\left((\pm (k-k_i)+k+k_i)T_1^i+dT_2^i\right)\\
&=dlF+\sum_{i=1}^{s}\left((\pm (k-k_i)+k+k_i)T_1^i+dT_2^i\right)\\
&\quad +\sum_{i=s+1}^{N}\left((\pm (k-2k_i)+k+2k_i)T_1^i+dT_2^i\right).
\end{split}
\end{equation}
Let $q$ be an element in $Q$. Note that this time
\begin{align*}
(k-k_i)+(k_i+k)&=d\\
-(k-k_i)+(k_i+k)&=d_i
\end{align*}
for $i\leq s$ while 
\begin{align*}
(k-2k_i)+(k+2k_i)&=d\\
-(k-2k_i)+(k+2k_i)&=2d_i
\end{align*}
for $i\geq s+1$. Since $q$ is the greatest common divisor of certain elements $d_i$ for $i\leq s$ and $2d_i$ for $i\geq s+1$, it follows as above that we can choose the plus and minus signs appropriately to get a symplectic structure $\omega_I$ on $W$ whose canonical class has divisibility equal to $q$.   
\end{proof}

\begin{ex} Suppose that $d=45$ and choose $d_0=45, d_1=15, d_2=9, d_3=5$. Then $Q=\{45,15,9,5,3,1\}$ and for every integer $n\geq 7$ there exists a homotopy elliptic surfaces $W$ with $\chi_h(W)=n$ that admits at least $6$ inequivalent symplectic structures whose canonical classes have divisibilities given by the elements in $Q$. One can also find an infinite family of homeomorphic but non-diffeomorphic manifolds of this kind. 
\end{ex}
\begin{cor}\label{cor 2 to the m symp on homotp ellipt} Let $m\geq 1$ be an arbitrary integer.
\begin{itemize}
\item There exist simply-connected non-spin 4-manifolds $W$ homeomorphic to the elliptic surfaces $E(2m+1)$ and $E(2m+2)_2$ which admit at least $2^m$ inequivalent symplectic structures. 
\item There exist simply-connected spin 4-manifolds $W$ homeomorphic to $E(6m-2)$ and $E(6m)$ which admit at least $2^{2m-1}$ inequivalent symplectic structures and spin manifolds homeomorphic to $E(6m+2)$ which admit at least $2^{2m}$ inequivalent symplectic structures.
\end{itemize}
\end{cor}
\begin{proof} Choose $N$ pairwise different odd prime numbers $p_1,\dotsc,p_N$. Let $d=d_0=p_1\cdot\dotsc\cdot p_N$ and consider the integers
\begin{align*}
d_1&=p_2\cdot p_3\cdot \dotsc\cdot p_N\\
d_2&=p_1\cdot p_3\cdot \dotsc\cdot p_N\\
&\vdots\\
d_N&=p_1\cdot \dotsc \cdot p_{N-1},
\end{align*}
obtained by deleting the corresponding prime in $d$. Then the associated set $Q$ of greatest common divisors consists of all products of the $p_i$ where each prime occurs at most once: If such a product $x$ does not contain precisely the primes $p_{i_1},\dotsc,p_{i_r}$ then $x$ is the greatest common divisor of $d_{i_1},\dotsc,d_{i_r}$. The set $Q$ has $2^N$ elements. 

Let $m\geq 1$ be an arbitary integer. Setting $N=m$ there exists by Theorem \ref{elliptic real q} for every integer $n\geq 2N+1=2m+1$ a homotopy elliptic surface $W$ with $\chi_h(W)=n$ which has $2^m$ symplectic structures realizing all elements in $Q$ as the divisibility of their canonical classes. Since $d$ is odd, the 4-manifolds $W$ are non-spin.

Setting $N=2m-1$ there exists for every even integer $n\geq 3N+1=6m-2$ a homotopy elliptic surface $W$ with $\chi_h(W)=n$ which has $2^{2m-1}$ symplectic structures realizing all elements in $Q$ multiplied by $2$ as the divisibility of their canonical classes. Since all divisibilities are even, the manifold $W$ is spin. Setting $N=2m$ we can choose $n=6m+2$ to get a spin homotopy elliptic surface $W$ with $\chi_h(W)=6m+2$ and $2^{2m}$ inequivalent symplectic structures.
\end{proof} 

The construction in the proof of Theorem \ref{elliptic real q} can be extended to the spin manifolds in Theorem \ref{thm existence c1 positive d even} with $c_1^2>0$: 
\begin{thm}\label{real q c12 positive d even} Let $N\geq 1$ be an integer. Suppose that $d\geq 2$ is an even integer and $d_0,\dotsc,d_N$ are positive even integers dividing $d$ as in Definition \ref{defn set Q of divisors}. Let $Q$ be the associated set of greatest common divisors. Let $m$ be an integer such that $2m\geq 3N+2$ and $t\geq 1$ an arbitrary integer. Then there exists a simply-connected closed spin 4-manifold $W$ with invariants
\begin{align*}
c_1^2(W)&=2td^2\\
e(W)&=td^2+24m\\
\sigma(W)&=-16m,
\end{align*}
and the following property: For each integer $q\in Q$ the manifold $W$ admits a symplectic structure whose canonical class $K$ has divisibility equal to $q$. Hence $W$ admits at least $|Q|$ many inequivalent symplectic structures.
\end{thm}  
\begin{proof} Let $l=2m-2N$. By the construction of Theorem \ref{thm existence c1 positive d even} there exists a simply-connected symplectic spin 4-manifold $X$ with invariants
\begin{align*}
c_1^2(X)&= 2td^2\\
e(X)&= td^2+12l\\
\sigma(X)&= -8l\\
K_X&= d(\tfrac{1}{2}lF+R+t\Sigma_M).
\end{align*}
In particular, the canonical class of $X$ has divisibility $d$. In the construction of $X$ starting from the elliptic surface $E(l)$ we have only used one Lagrangian rim torus. Hence $l-2$ of the $l-1$ triples of Lagrangian rim tori in $E(l)$ (cf.~Example \ref{n-1 Lagrangian triples in E(n)}) remain unchanged. Note that $l-2\geq N$ by our assumptions. Since the symplectic form on $E(l)$ in a neighbourhood of these tori does not change in the construction of $X$ by the Gompf fibre sum, we can assume that $X$ contains at least $N$ triples of Lagrangian tori as in the proof of Theorem \ref{elliptic real q}. We can now use the same construction as in this theorem on the $N$ triples of Lagrangian tori in $X$ to get a simply-connected spin 4-manifold $W$ with invariants
\begin{align*}
c_1^2(W)&= 2td^2\\
e(W)&= td^2+12l+24N=td^2+24m\\
\sigma(W)&= -8l-16N=-16m,\\
\end{align*}
admitting $2^N$ symplectic structures. In particular, for each $q\in Q$ the manifold $W$ admits a symplectic structure $\omega_I$ whose canonical class is given by the formulas in equation \eqref{canonical class of W pos c12 d not div 4} and \eqref{canonical class of W pos c12 d does div 4} where the term $dlF$ is replaced by $K_X=d(\tfrac{1}{2}lF+R+t\Sigma_M)$. It follows again that the canonical class of $\omega_I$ has divisibility precisely equal to $q$. 
\end{proof}

\begin{cor}\label{cor two symp inequiv spin c12 pos} Let $d\geq 6$ be an even integer and $t\geq 1, m\geq 3$ arbitrary integers. Then there exists a simply-connected closed spin 4-manifold $W$ with invariants
\begin{align*}
c_1^2(W)&=2td^2\\
e(W)&=td^2+24m\\
\sigma(W)&=-16m,
\end{align*}
such that $W$ admits at least two inequivalent symplectic structures. 
\end{cor}
This follows with $N=1$ and choosing $d_0=d$ and $d_1=2$, since in this case $Q$ consists of two elements.
\begin{ex}\label{horikawa example several symp} We consider Corollary \ref{cor two symp inequiv spin c12 pos} for the spin homotopy Horikawa surfaces in Example \ref{horikawa with div K}. Let $t\geq 1$ and $k\geq 3$ be arbitrary odd integers and define an integer $m$ by $2m=3tk^2+3$. Let $d=2k$ and $d_1=2$. Since $d=2k$ is not divisible by $4$, the set $Q$ is equal to $\{2k, 2\}$ by Definition \ref{defn set Q of divisors}. Hence there exists a spin homotopy Horikawa surface $X$ on the Noether line with invariants
\begin{align*}
c_1^2(X)&=8tk^2\\
\chi_h(X)&=4tk^2+3,
\end{align*}
which admits two inequivalent symplectic structures: the canonical class of the first symplectic structure has divisibility $2k$ while the canonical class of the second symplectic structure is divisible only by $2$. 
\end{ex}

Similarly the construction in Theorem \ref{real q c12 positive d even} can be extended to the non-spin manifolds in Theorem \ref{odd div sign div 8} with $c_1^2>0$: 
\begin{thm}\label{real q c12 positive d odd} Let $N \geq 1$ be an integer. Suppose that $d\geq 3$ is an odd integer and $d_0,\dotsc,d_N$ positive integers dividing $d$ as in Definition \ref{defn set Q of divisors}. Let $Q$ be the associated set of greatest common divisors. Let $m\geq 2N+2$ and $t\geq 1$ be arbitrary integers. Then there exists a simply-connected closed non-spin 4-manifold $W$ with invariants
\begin{align*}
c_1^2(W)&=8td^2\\
e(W)&=4td^2+12m\\
\sigma(W)&=-8m,
\end{align*}
and the following property: For each integer $q\in Q$ the manifold $W$ admits a symplectic structure whose canonical class $K$ has divisibility equal to $q$. Hence $W$ admits at least $|Q|$ many inequivalent symplectic structures.
\end{thm}  
\begin{proof} The proof is analogous to the proof of Theorem \ref{real q c12 positive d even}. Let $l=m-N$. By the construction of Theorem \ref{odd div sign div 8} there exists a simply-connected non-spin symplectic 4-manifold $X$ with invariants
\begin{align*}
c_1^2(X)&=8td^2\\
e(X)&=4td^2+12l\\
\sigma(X)&=-8l,\\
\end{align*}
whose canonical class $K_X$ has divisibility $d$. The manifold $X$ contains $l-2$ triples of Lagrangian tori. By our assumptions $l-2\geq N$. Hence we can do the construction in Theorem \ref{elliptic real q} (for $d$ odd) to get a simply-connected non-spin 4-manifold $W$ with invariants
\begin{align*}
c_1^2(W)&=8td^2\\
e(W)&=4td^2+12l+12N=4td^2+12m\\
\sigma(X)&=-8l-8N=-8m.\\
\end{align*}
The 4-manifold $W$ admits for every integer $q\in Q$ a symplectic structure whose canonical class has divisibility equal to $q$.
\end{proof}
Choosing $N=1$, $d_0=d$ and $d_1=1$ the set $Q$ contains two elements. This implies:
\begin{cor} Let $d\geq 3$ be an odd integer and $t\geq 1, m\geq 4$ arbitrary integers. Then there exists a simply-connected closed non-spin 4-manifold $W$ with invariants
\begin{align*}
c_1^2(W)&=8td^2\\
e(W)&=4td^2+12m\\
\sigma(W)&=-8m,
\end{align*}
such that $W$ admits at least two inequivalent symplectic structures. 
\end{cor}

\section{Branched coverings}\label{sect branched coverings}

Let $M^n$ be a closed, oriented smooth manifold and $F^{n-2}$ a closed, oriented submanifold of codimension 2. Suppose that the fundamental class $[F]\in H_{n-2}(M;\mathbb{Z})$ is divisible by an integer $m>1$ and choose a class $B\in H_{n-2}(M;\mathbb{Z})$ such that $[F]=mB$. The integer $m$ together with the class $B$ determine a branched covering of $M$. 
\begin{defn} We denote by $\phi:M(F,B,m)\rightarrow M$ the $m$-fold branched covering of $M$ branched over $F$ and determined by $m$ and $B$.
\end{defn}
For the construction of branched coverings see Hirzebruch's article \cite{Hirz}. The smooth manifold $M(F,B,m)$ has the following properties:
\begin{itemize}
\item Over the complement $M'=M\setminus F$, the map $\phi:\phi^{-1}(M')\rightarrow M'$ is a standard $m$-fold cyclic covering. 
\item $\phi$ maps the submanifold $\overline{F}=\phi^{-1}(F)$ diffeomorphically onto $F$ and on tubular neighbourhoods, $\phi:\nu(\overline{F})\rightarrow \nu(F)$ is locally of the form
\begin{equation*}
U\times D^2\rightarrow U\times D^2, (x,z)\mapsto (x,z^m),
\end{equation*}
where $D^2$ is considered as the unit disk in $\mathbb{C}$. 
\end{itemize}
Suppose that $M$ is a smooth complex algebraic surface and $D\subset M$ a smooth connected complex curve. If $m>0$ is an integer that divides $[D]$ and $B\in H_2(M;\mathbb{Z})$ a homology class such that $[D]=mB$, then the branched covering $M(D,B,m)$ also admits the structure of an algebraic surface. Some invariants of this surface are given by the following proposition.  
\begin{prop}\label{inv branched cov M,D,m} Let $D$ be a smooth connected complex curve in a complex surface $M$ such that $[D]=mB$. Let $\phi\colon M(D,B,m)\rightarrow M$ be the branched covering. Then the invariants of $N:=M(D,B,m)$ are given by:   
\begin{enumerate}
\item $K_N=\phi^*(K_M+(m-1)B)$ 
\item $c_1^2(N) = m(K_M+(m-1)B)^2$ 
\item $e(N) = me(M)-(m-1)e(D)$,
\end{enumerate}
where $e(D)=2-2g(D)=-(K_M\cdot D+D^2)$ by the adjunction formula. 
\end{prop}
\begin{proof} The formula for $e(N)$ follows by a well-known formula for the Euler characteristic of a topological space decomposed into two pieces and the formula for standard, unramified coverings. The formula for $c_1^2(N)$ then follows by the signature formula of Hirzebruch \cite{Hirz}:
\begin{equation*}
\sigma(N)=m\sigma(M)-\frac{m^2-1}{3m}D^2.
\end{equation*}
Finally, the formula for $K_N$ can be found in \cite[Chapter I, Lemma 17.1]{BPV}. 
\end{proof}

Suppose that the complex curve $D$ is contained in the linear system $|nK_M|$ and hence represents in homology a multiple $nK_M$ of the canonical class of $M$. Let $m>0$ be an integer dividing $n$ and write $n=ma$. Setting $[D]=nK_M$ and $B=aK_M$ in Proposition \ref{inv branched cov M,D,m} we get.
\begin{cor} \label{branched cov inv} Let $D$ be a smooth connected complex curve in a complex surface $M$ with $[D]=nK_M$ and $\phi\colon M(D,aK_M,m)\rightarrow M$ the branched covering. Then the invariants of $N:=M(D,aK_M,m)$ are given by:
\begin{enumerate}
\item $K_N = (n+1-a)\phi^*K_M$
\item $c_1^2(N)=m(n+1-a)^2c_1^2(M)$
\item $e(N)=me(M)+(m-1)n(n+1)c_1^2(M)$
\end{enumerate}
\end{cor}
 
We consider again the general situation that $M$ is a smooth, oriented manifold and $F$ an oriented submanifold of codimension $2$. The fundamental group of $M$ is related to the fundamental group of the complement $M'=M\setminus F$ by
\begin{equation}\label{fundament group complement F}
\pi_1(M)\cong \pi_1(M')/N(\sigma),
\end{equation}
where $\sigma$ denotes the meridian to $F$, given by a circle fibre of $\partial \nu(F)\rightarrow F$, and $N(\sigma)$ denotes the normal subgroup in $\pi_1(M')$ generated by this element (a proof can be found in the appendix of \cite{MHthesis}). Using this formula the fundamental group of a branched covering can be calculated in the following case.
\begin{thm}\label{fund branch cov} Let $M^n$ be a closed oriented manifold and $F^{n-2}$ a closed oriented submanifold. Suppose in addition that the fundamental group of $M'$ is abelian. Then for all $m$ and $B$ with $[F]=mB$ there exists an isomorphism
\begin{equation*}
\pi_1(M(F,B,m))\cong\pi_1(M).
\end{equation*}
\end{thm}
\begin{proof} Let $k>0$ denote the maximal integer dividing $[F]$. Since $m$ divides $k$, we can write $k=ma$ with $a>0$. Let $\overline{M'}$ denote the complement to $\overline{F}$ in $M(F,B,m)$ and $\overline{\sigma}$ the meridian to $\overline{F}$. By equation \eqref{fundament group complement F} we have
\begin{equation*}
\pi_1(M(F,B,m))\cong \pi_1(\overline{M'})/N(\overline{\sigma}).
\end{equation*}
There is an exact sequence
\begin{equation*}
0\rightarrow \pi_1(\overline{M'})\stackrel{\pi_*}{\rightarrow} \pi_1(M')\rightarrow \mathbb{Z}_m\rightarrow 0,
\end{equation*}
since $\pi\colon \overline{M'}\rightarrow M'$ is an $m$-fold cyclic covering. The assumption that $\pi_1(M')$ is abelian implies that $\pi_1(\overline{M'})$ is also abelian. Therefore the normal subgroups generated by the meridians in these groups are cyclic and we get an exact sequence of subgroups
\begin{equation*}
0\rightarrow \mathbb{Z}_a\overline{\sigma}\stackrel{m\cdot}{\rightarrow}\mathbb{Z}_{ma}\sigma\rightarrow \mathbb{Z}_m\rightarrow 0,
\end{equation*}
where $\sigma$ is the meridian to $F$. The third map being a surjection implies that for every element $\alpha\in \pi_1(M')$ there is an integer $r\in\mathbb{Z}$ such that $\alpha+r\sigma$ maps to zero in $\mathbb{Z}_m$ and hence is in the image of $\pi_*$. In other words, the induced map
\begin{equation*}
\pi_*\colon \pi_1(\overline{M'})\longrightarrow \pi_1(M')/\langle\sigma\rangle
\end{equation*}
is surjective. The kernel of this map is equal to $\langle\overline{\sigma}\rangle$, hence
\begin{equation*}
\pi_1(\overline{M'})/\langle\overline{\sigma}\rangle\stackrel{\cong}{\longrightarrow}\pi_1(M')/\langle\sigma\rangle.
\end{equation*}
Again by equation \eqref{fundament group complement F} this implies $\pi_1(M(F,B,m))\cong\pi_1(M)$. 
\end{proof}

We want to apply this theorem in the case where $M$ is a 4-manifold and $F$ is an embedded surface. Even if $M$ is simply-connected the complement $M'$ does not have abelian fundamental group in general. However, in the complex case, we can use the following theorem of Nori \cite[Proposition 3.27]{N}.
\begin{thm}\label{nori thm} Let $M$ be a smooth complex algebraic surface and $D,E\subset M$ smooth complex curves which intersect transversely. Assume that $D'^2>0$ for every connected component $D'\subset D$. Then the kernel of $\pi_1(M\setminus (D\cup E))\rightarrow \pi_1(M\setminus E)$ is a finitely generated abelian group.
\end{thm}
In particular for $E=\emptyset$, this implies that the kernel of
\begin{equation*}
\pi_1(M')\rightarrow \pi_1(M)
\end{equation*}
is a finitely generated abelian group if $D$ is connected and $D^2>0$, where $M'$ denotes $M\setminus D$. If $M$ is simply-connected, it follows that $\pi_1(M')$ is abelian. Hence with Theorem \ref{fund branch cov} we get:  
 
\begin{cor}\label{simple branched} Let $M$ be a simply-connected, smooth complex algebraic surface and $D\subset M$ a smooth connected complex curve with $D^2>0$. Let $\overline{M}$ be a cyclic ramified cover of $M$ branched over $D$. Then $\overline{M}$ is also simply-connected.  
\end{cor}
Catanese \cite{Cat2} has also used in a different situation restrictions on divisors to ensure that certain ramified coverings are simply-connected.

\section{Surfaces of general type and pluricanonical systems}\label{sec geo surf gen}

In this section we collect some results concerning the geography of simply-connected surfaces of general type and the existence of smooth divisors in pluricanonical systems. 

The following result due to Persson \cite[Proposition 3.23]{Pe} is the main geography result we use for our constructions. 
\begin{thm}\label{geog Per I} Let $x,y$ be positive integers such that
\begin{equation*}
2x-6\leq y\leq 4x-8.
\end{equation*}
Then there exists a simply-connected minimal complex surface $M$ of general type such that $\chi_h(M)=x$ and $c_1^2(M)=y$. Furthermore, $M$ can be chosen as a genus 2 fibration. 
\end{thm}
The smallest integer $x$ to get an inequality which can be realized with $y>0$ is $x=3$. Since $\chi_h(M)=p_g(M)+1$ for simply-connected surfaces, this corresponds to surfaces with $p_g=2$. Hence from Theorem \ref{geog Per I} we get minimal simply-connected complex surfaces $M$ with
\begin{equation*}
p_g=2\,\,\mbox{and}\,K^2=1,2,3,4.
\end{equation*}
Similarly for $x=4$ we get surfaces with
\begin{equation*}
p_g=3\,\,\mbox{and}\,K^2=2,\dotsc,8.
\end{equation*}
For surfaces with $K^2=1$ or $K^2=2$ we have:
\begin{prop}\label{small K^2 realized} For $K^2=1$ and $K^2=2$ all possible values for $p_g$ given by the Noether inequality $K^2\geq 2p_g-4$ can be realized by simply-connected minimal complex surfaces of general type.
\end{prop}
\begin{proof}
By the Noether inequality, only the following values for $p_g$ are possible:
\begin{align*}
K^2=1:&\quad p_g=0,1,2\\
K^2=2:&\quad p_g=0,1,2,3.
\end{align*}
The cases $K^2=1, p_g=2$ and $K^2=2, p_g=2,3$ are covered by Persson's theorem. In particular, the surfaces with $K^2=1, p_g=2$ and $K^2=2, p_g=3$ are Horikawa surfaces described in \cite{HorI, HorII}. The remaining cases can also be covered: The Barlow surface, constructed in \cite{Barl}, is a simply-connected numerical Godeaux surface, i.e.~a minimal complex surface of general type with $K^2=1, p_g=0$. Simply-connected minimal surfaces of general type with $K^2=1,2$ and $p_g=1$ exist by constructions due to Enriques, see \cite{Cat, CP, Cha}. Finally, Lee and Park have constructed in \cite{LP} a simply-connected minimal surface of general type with $K^2=2, p_g=0$. It is a numerical Campedelli surface. 
\end{proof}
Suppose $M$ is a minimal smooth complex algebraic surface of general type and consider the multiples $L=nK=K^{\otimes n}$ of the canonical line bundle of $M$. By a theorem of Bombieri \cite{Bom2, BPV} all divisors in the linear system $|nK|$ are connected. If $|nK|$ has no fixed parts and is base point free, it determines an everywhere defined holomorphic map to a projective space and we can find a non-singular divisor representing $nK$ by taking the preimage of a generic hyperplane section.
\begin{thm}\label{pluri real} Let $M$ be a minimal smooth complex algebraic surface of general type. Then the pluricanonical system $|nK|$ determines an everywhere defined holomorphic map in the following cases:
\begin{itemize}
\item $n\geq 4$
\item $n=3$ and $K^2\geq 2$
\item $n=2$ and $K^2\geq 5$ or $p_g\geq 1$.
\end{itemize}
\end{thm}
For proofs and references see \cite{Bom, Bom2, CT, Kod, MP, Re}. 
\begin{rem}\label{rem godeaux 3K base point} In some of the remaining cases it is also known that pluricanonical systems define a holomorphic map. In particular, suppose that $M$ is a numerical Godeaux surface. Then the map defined by $|3K|$ is holomorphic if $H_1(M;\mathbb{Z})=0$ or $\mathbb{Z}_2$, e.g.\ if $M$ is simply-connected \cite{Miy, MRe}. This is also known for the map defined by $|2K|$ in the case of a simply-connected surface $M$ with $K^2=4, p_g=0$ by \cite{CT, Ko1}.
\end{rem}

\section{Branched covering construction of algebraic surfaces with divisible canonical class}\label{section branched cov with div}

Suppose that $M$ is a simply-connected minimal complex surface of general type. Let $m,d\geq 2$ be integers such that $m-1$ divides $d-1$ and define the integers
\begin{equation*}
a=\tfrac{d-1}{m-1}\quad\mbox{and}\quad n=ma.
\end{equation*}
Then $d=n+1-a$ and the assumptions imply that $n\geq 2$. We assume in addition that $nK_M$ can be represented by a smooth complex connected curve $D$ in $M$, see Theorem \ref{pluri real}. Let $\overline{M}=M(D, aK_M,m)$ denote the associated $m$-fold branched cover over the curve $D$. 
\begin{thm}\label{invariants of plurican branched prop} Let $M$ be a simply-connected minimal surface of general type and $m,d\geq 2$ integers such that $d-1$ is divisible by $m-1$ with quotient $a$. Suppose that $D$ is a smooth connected curve in the linear system $|nK_M|$ where $n=ma$. Then the $m$-fold cover of $M$, branched over $D$, is a simply-connected complex surface $\overline{M}$ of general type with invariants
\begin{itemize}
\item $K_{\overline{M}} = d\phi^*K_M$
\item $c_1^2(\overline{M})=md^2c_1^2(M)$
\item $e(\overline{M})=m(e(M)+(d-1)(d+a)c_1^2(M))$
\item $\chi_h(\overline{M})=m\chi_h(M)+\frac{1}{12}m(d-1)(2d+a+1)c_1^2(M)$
\item $\sigma(\overline{M})=-\frac{1}{3}m(2e(M)+(d(d-2)+2a(d-1))c_1^2(M))$.
\end{itemize}
In particular, the canonical class $K_{\overline{M}}$ is divisible by $d$ and $\overline{M}$ is minimal. 
\end{thm}
\begin{proof}
The formulas for the invariants are given by Corollary \ref{branched cov inv}. Since 
$D^2=n^2K_M^2>0$, the complex surface $\overline{M}$ is simply-connected by Corollary \ref{simple branched}. Moreover, $\overline{M}$ is of general type because $c_1^2(\overline{M})>0$ and $\overline{M}$ cannot be rational or ruled. The claim about minimality follows from Lemma \ref{lem 2 geography K}, because the divisibility of $K_{\overline{M}}$ is at least $d\geq 2$.
\end{proof}
Note that the signature $\sigma(\overline{M})$ is always negative, hence surfaces with positive signature cannot be constructed in this way. 

The transformation 
\begin{equation*}
\Phi\colon (e(M),c_1^2(M))\mapsto (e(\overline{M}), c_1^2(\overline{M}))
\end{equation*}
given by Theorem \ref{invariants of plurican branched prop} is linear and can be written as
\begin{equation*}
\left(\begin{array}{c} e(\overline{M}) \\ c_1^2(\overline{M}) \end{array}\right)=m\left(\begin{array}{cc} 1&\Delta \\ 0&d^2 \end{array}\right)\left(\begin{array}{c} e(M) \\ c_1^2(M) \end{array}\right),
\end{equation*}
with the abbreviation 
\begin{equation*}
\Delta=(d-1)(d+a).
\end{equation*}
This map is invertible over $\mathbb{R}$ and maps the quadrant $\mathbb{R}^+\times\mathbb{R}^+$ into the same quadrant, where $\mathbb{R}^+\times\mathbb{R}^+$ denotes the quadrant in the $\mathbb{R}\times\mathbb{R}$--plane where both coordinates have positive entries. The inverse of $\Phi$ is given by
\begin{equation*}
\left(\begin{array}{c} e(M) \\ c_1^2(M) \end{array}\right)=\tfrac{1}{m}\left(\begin{array}{cc} 1&-\Delta/d^2 \\ 0&1/d^2\end{array}\right)\left(\begin{array}{c} e(\overline{M}) \\ c_1^2(\overline{M}) \end{array}\right).
\end{equation*}
\begin{defn} We call a point in the $\mathbb{R}^+\times\mathbb{R}^+$ quadrant {\em admissible} if $e(M)+c_1^2(M)\equiv 0$ mod $12$.
\end{defn}
The coordinates $e(M)$ and $c_1^2(M)$ of a complex surface are always admissible by the Noether formula.
\begin{lem}\label{lem image Phi admissible points} The image of the admissible points in the quadrant $\mathbb{R}^+\times\mathbb{R}^+$ under the map $\Phi$ consists of the points satisfying $e(\overline{M})\equiv 0$ mod $m$, $c_1^2(\overline{M})\equiv 0$ mod $md^2$ and 
\begin{equation*}
\mbox{$\tfrac{1}{m}e(\overline{M})+\tfrac{1-\Delta}{md^2}c_1^2(\overline{M})\equiv 0$ mod $12$.}
\end{equation*}
\end{lem}
The proof is immediate by the formula for the inverse of $\Phi$. We want to calculate the image under $\Phi$ of the sector given by Theorem \ref{geog Per I}. First, we rewrite Persson's theorem in the following equivalent form (we omit the proof):
\begin{cor}\label{cor equiv persson sector} Let $e, c$ be positive integers such that $c\geq 36-e$ and $e+c\equiv 0$ mod $12$. If
\begin{equation*}
\tfrac{1}{5}(e-36)\leq c\leq \tfrac{1}{2}(e-24),
\end{equation*}
then there exists a simply-connected minimal surface $M$ of general type with invariants $e(M)=e$ and $c_1^2(M)=c$.
\end{cor}
In the next step, we calculate the image under $\Phi$ of the lines in the $(e,c)$--plane which appear in this corollary. A short calculation shows that the line $c=\frac{1}{5}(e-36)$ maps to
\begin{equation}\label{new 2x-6 line}
c_1^2(\overline{M})= \frac{d^2}{(5+\Delta)}(e(\overline{M})-36m),
\end{equation}
while the line $c=\frac{1}{2}(e-24)$ maps to
\begin{equation}\label{new 4x-8 line}
c_1^2(\overline{M})=\frac{d^2}{(2+\Delta)}(e(\overline{M})-24m).
\end{equation}
Similarly, the constraint $c\geq 36-e$ maps to
\begin{equation}\label{new 3 line}
c_1^2(\overline{M})\geq-\frac{d^2}{(1-\Delta)}(e(\overline{M})-36m).
\end{equation}
It follows that the image under $\Phi$ of the lattice points given by the constraints in Corollary \ref{cor equiv persson sector} consists precisely of the points in the sector between the lines \eqref{new 2x-6 line} and \eqref{new 4x-8 line}, which satisfy the constraint \eqref{new 3 line} and the constraints in Lemma \ref{lem image Phi admissible points}.

The surfaces in Theorem \ref{geog Per I} satisfy $p_g\geq 2$ and $K^2\geq 1$. By Theorem \ref{pluri real}, the linear system $|nK|$ for $n\geq 2$ on these surfaces defines a holomorphic map, except possibly in the case $p_g=2, K^2=1$ and $n=3$. Since $n=ma$ and $m\geq 2$, this occurs only for $m=3, a=1$ and $d=3$. The corresponding image under $\Phi$ has invariants $(e,c_1^2)=(129,27)$. This exception is implicitly understood in the following theorem. In all other cases we can consider the branched covering construction above. This can be summarized as follows: Consider integers $m,a,d$ as above, with $m,d\geq 2, a\geq 1$ and $\Delta=(d-1)(d+a)$. 
\begin{thm}\label{thm image persson branched} Let $x,y$ be positive integers such that $y(1-\Delta)\geq 36-x$ and $x+(1-\Delta)y\equiv 0$ mod $12$. If 
\begin{equation*}
\frac{1}{(5+\Delta)}(x-36)\leq y\leq \frac{1}{(2+\Delta)}(x-24),
\end{equation*} 
then there exists a simply-connected minimal complex surface $\overline{M}$ of general type with invariants $e(\overline{M})=mx$ and $c_1^2(\overline{M})=md^2y$, such that the canonical class of $\overline{M}$ is divisible by $d$. 
\end{thm}
We calculate some explicit examples for the branched covering construction given by Theorem \ref{thm image persson branched} and for some additional surfaces not covered by Persson's theorem. For any $d\geq 2$, we can choose $m=2$ and $a=d-1$ corresponding to 2-fold covers branched over $(2d-2)K$. The formulas for the invariants simplify to
\begin{itemize}
\item $c_1^2(\overline{M})=2d^2c_1^2(M)$
\item $e(\overline{M})=24\chi_h(M)+2d(2d-3)c_1^2(M)$
\item $\chi_h(\overline{M})=2\chi_h(M)+\frac{1}{2}d(d-1)c_1^2(M)$.
\end{itemize}
The first two examples are double coverings with $m=2$, the third example uses coverings of higher degree. Because of their topological invariants some of the surfaces are homeomorphic by Freedman's theorem to simply-connected symplectic 4-manifolds constructed in Sections \ref{sect geo spin neg sigma c pos} and \ref{sect geo non-spin neg sigma c pos}.

\begin{ex} We consider the Horikawa surfaces \cite{HorI} on the Noether line $c_1^2=2\chi_h-6$, which exist for every $\chi_h\geq 4$ and are also given by Persson's theorem \ref{geog Per I}. In this case $p_g\geq 3$ and $c_1^2\geq 2$, hence by Theorem \ref{pluri real} the linear system $|nK|$ for $n\geq 2$ defines a holomorphic map on these surfaces.
\begin{prop} Let $M$ be a Horikawa surface on the Noether line $c_1^2=2\chi_h-6$ where $\chi_h=4+l$ for $l\geq 0$. Then the 2-fold cover $\overline{M}$ of the surface $M$, branched over $(2d-2)K_M$ for an integer $d\geq 2$, has invariants
\begin{itemize}
\item $c_1^2(\overline{M})=4d^2(l+1)$
\item $\chi_h(\overline{M})=6+(2+d(d-1))(l+1)$
\item $e(\overline{M})=72+4(l+1)(6+2d^2-3d)$
\item $\sigma(\overline{M})=-48-4(l+1)(4+d^2-2d)$.
\end{itemize}
The canonical class $K_{\overline{M}}$ is divisible by $d$.
\end{prop}
For $d$ even, the integer $d^2-2d=d(d-2)$ is divisible by $4$, hence $\sigma$ is indeed divisible by $16$, which is necessary by Rochlin's theorem. The invariants are on the line
\begin{equation*}
c_1^2(\overline{M})=\frac{4d^2}{2+d(d-1)}(\chi_h(\overline{M})-6),
\end{equation*}
which has inclination close to $4$ for $d$ very large. 
\end{ex}

\begin{ex}\label{small c_1^2 div alg surfaces} We calculate the invariants for the branched covers with $m=2$ and integers $d\geq 3$ for the surfaces given by Proposition \ref{small K^2 realized}. Since $n=ma\geq 4$ in this case, Theorem \ref{pluri real} shows that the linear system $|nK|$ defines a holomorphic map and we can use the branched covering construction. 
\begin{prop} Let $M$ be a minimal complex surface of general type with $K^2=1$ or $2$. Then the 2-fold cover $\overline{M}$ of the surface $M$, branched over $(2d-2)K_M$ for an integer $d\geq 3$, has invariants
\begin{enumerate}
\item If $K^2=1$ and $p_g=0,1,2:$
\begin{itemize} 
\item $c_1^2(\overline{M})=2d^2$
\item $e(\overline{M})=24(p_g+1)+2d(2d-3)$
\item $\sigma(\overline{M})=-16(p_g+1)-2d(d-2)$.
\end{itemize}
\item If $K^2=2$ and $p_g=0,1,2,3:$
\begin{itemize} 
\item $c_1^2(\overline{M})=4d^2$
\item $e(\overline{M})=24(p_g+1)+4d(2d-3)$
\item $\sigma(\overline{M})=-16(p_g+1)-4d(d-2)$.
\end{itemize}
\end{enumerate} 
In both cases the canonical class $K_{\overline{M}}$ is divisible by $d$.
\end{prop}
\end{ex}

\begin{ex} Consider the Barlow surface $M_B$ and the surface $M_{LP}$ of Lee and Park, mentioned in the proof of Proposition \ref{small K^2 realized}. The invariants are
\begin{itemize}
\item $c_1^2(M_B)=1, \chi_h(M_B)=1$ and $e(M_B)=11$
\item $c_1^2(M_{LP})=2,\chi_h(M_{LP})=1$ and $e(M_{LP})=10$.
\end{itemize}
By Theorem \ref{pluri real}, we can consider branched covers over both surfaces with $ma\geq 3$ (the Barlow surface is a simply-connected numerical Godeaux surface, hence $|3K|$ defines a holomorphic map by Remark \ref{rem godeaux 3K base point}). See Tables \ref{Branched cov Barlow} and \ref{Branched cov Lee Park} for a calculation of the invariants of $\overline{M}$ for small values of $d$ and $m$. There is a coincidence between the $4$-fold cover of the Barlow surface branched over $4K_M$ and the $2$-fold cover of the surface of Lee and Park branched over $6K_M$: Both have the same Chern invariants and the same divisibility $d=4$ of the canonical class. Hence the manifolds are homeomorphic and by Theorem \ref{thm geog intro Kahler}, both branched coverings have the same Seiberg-Witten invariants.
\begin{table}[hbt]
\setlength{\tabcolsep}{0.16cm}
\renewcommand{\arraystretch}{1.5}
\begin{tabular}{|c|c|c|c|c|c|c|c|c|} \hline
 $d$ & $m$ &$ma$ &$(d-1)(d+a)$ &$e(\overline{M_B})$& $c_1^2(\overline{M_B})$  &$\chi_h(\overline{M_B})$& $b_2^+(\overline{M_B})$& $\sigma(\overline{M_B})$ \\ 
\hline\hline
 $3$ & $2$& $4$& $10$ &$42$ & $18$ & $5$ & $9$ & $-22$ \\ \hline
 $3$ & $3$ &$3$& $8$ & $57$& $27$ & $7$  & $13$ & $-29$\\ \hline
 $4$ & $2$ &$6$& $21$ & $64$ & $32$ & $8$ & $15$ & $-32$\\ \hline
$4$ & $4$ &$4$& $15$ &$104$ &$64$ & $14$ & $27$ & $-48$\\ \hline
 $5$ & $2$ &$8$& $36$ &$94$ &$50$& $12$ & $23$ & $-46$\\ \hline
$5$ & $3$ &$6$& $28$ &$117$& $75$& $16$  &  $31$ & $-53$ \\ \hline
 $5$ & $5$ &$5$& $24$ &$175$ &$125$ & $25$   & $49$ &$-75$ \\ \hline
 $6$ & $2$ &$10$& $55$ &$132$& $72$& $17$ & $33$ & $-64$ \\ \hline
 $6$ & $6$ &$6$& $35$ &$276$& $216$ & $41$ & $81$ &$-112$ \\ \hline
\end{tabular}
\caption{Ramified coverings of the Barlow surface $M_B$ of degree $m$ branched over $maK$.}\label{Branched cov Barlow}
\end{table}

\begin{table}[hbt]
\setlength{\tabcolsep}{0.11cm}
\renewcommand{\arraystretch}{1.5}
\begin{tabular}{|c|c|c|c|c|c|c|c|c|} \hline
$d$ & $m$ &$ma$ &$(d-1)(d+a)$ &$e(\overline{M_{LP}})$& $c_1^2(\overline{M_{LP}})$ & $\chi_h(\overline{M_{LP}})$ & $b_2^+(\overline{M_{LP}})$& $\sigma(\overline{M_{LP}})$ \\ 
\hline\hline
$3$ & $2$ & $4$ & $10$ & $60$ & $36$ & $8$  & $15$ & $-28$ \\ \hline
$3$ & $3$ & $3$ & $8$ & $78$ & $54$ & $11$ & $21$ & $-34$\\ \hline
$4$ & $2$ & $6$ & $21$ & $104$ & $64$ & $14$ & $27$ & $-48$\\ \hline
 $4$ & $4$ & $4$ & $15$ & $160$ & $128$ & $24$  & $47$ & $-64$\\ \hline
$5$ & $2$ & $8$ & $36$ & $164$ & $100$ & $22$  &  $43$ & $-76$\\ \hline
 $5$ & $3$ & $6$ & $28$ & $198$ & $150$ & $29$   & $57$ & $-82$ \\ \hline
 $5$ & $5$ & $5$ & $24$ & $290$ & $250$ & $45$ &   $89$ & $-110$ \\ \hline
 $6$ & $2$ & $10$ & $55$ & $240$ & $144$ & $32$  &  $63$ & $-112$ \\ \hline
 $6$ & $6$ & $6$ & $35$ & $480$ & $432$ & $76$   & $151$ &  $-176$ \\ \hline
\end{tabular}
\caption{Ramified coverings of the Lee-Park surface $M_{LP}$ of degree $m$ branched over $maK$.}\label{Branched cov Lee Park}
\end{table}
\end{ex}

\begin{rem}\label{cov sing curve}
More general examples are possible by considering branched coverings over singular complex curves. The following example is described for instance in \cite[Chapter 7]{GS}: Let $B_{n,m}$ denote the singular complex curve in $\mathbb{C}P^1\times \mathbb{C}P^1$, which is the union of $2n$ parallel copies of the first factor and $2m$ parallel copies of the second factor. The curve $B_{n,m}$ represents in cohomology the class $2nS_1+2mS_2$, where $S_1=[\mathbb{C}P^1\times \{*\}]$ and $S_2=[\{*\}\times \mathbb{C}P^1]$. Let $X'(n,m)$ denote the double covering of $\mathbb{C}P^1\times \mathbb{C}P^1$ branched over $B_{n,m}$. It is a singular complex surface, which has a canonical resolution $X(n,m)$ (see \cite[Chapter III]{BPV}). As a smooth 4-manifold, $X(n,m)$ is diffeomorphic to the double cover of $\mathbb{C}P^1\times \mathbb{C}P^1$ branched over the smooth curve $\widetilde{B}_{n,m}$ given by smoothing the double points. Hence the topological invariants for $X=X(n,m)$ can be calculated with the formulas from Proposition \ref{inv branched cov M,D,m}:
\begin{itemize}
\item $c_1^2(X)=4(n-2)(m-2)$
\item $e(X)=6+2(2m-1)(2n-1)$
\item $\sigma(X)=-4mn$
\end{itemize}
Writing $X'=X'(n,m)$ and $M=\mathbb{C}P^1\times \mathbb{C}P^1$, let $\phi\colon X'\rightarrow M$ denote the double covering, $\pi\colon X\rightarrow X'$ the canonical resolution and $\psi=\phi\circ \pi$ the composition. Since all singularities of $B_{n,m}$ are ordinary double points, the canonical class of $X$ can be calculated by a formula in \cite[Theorem 7.2, Chapter III]{BPV}:
\begin{align*}
K_X&=\psi^*(K_M+\tfrac{1}{2}B_{m,n})\\
&= \psi^*(-2S_1-2S_2+nS_1+mS_2)\\
&= \psi^*((n-2)S_1+(m-2)S_2).
\end{align*}
This formula has the following interpretation: The map $\psi\colon X\rightarrow \mathbb{C}P^1\times \mathbb{C}P^1$ followed by the projection onto the first factor defines a fibration $X\rightarrow \mathbb{C}P^1$ whose fibres are the branched covers of the rational curves $\{p\}\times \mathbb{C}P^1$, where $p\in \mathbb{C}P^1$. The generic rational curve among them is disjoint from the $2m$ curves in $B_{n,m}$ parallel to $\{*\}\times \mathbb{C}P^1$ and intersects the $2n$ curves parallel to $\mathbb{C}P^1\times \{*\}$ in $2n$ points. This implies that the generic fibre $F_2$ of the fibration is a double branched cover of $\mathbb{C}P^1$ in $2n$ distinct points and hence a smooth complex curve of genus $n-1$. This curve represents the class $\psi^*S_2$ in the surface $X$. Similarly, there is a fibration $X\rightarrow \mathbb{C}P^1$ in genus $m-1$ curves which represent $F_1=\psi^*S_1$. Hence we can write
\begin{equation*}
K_X=(n-2)F_1+(m-2)F_2.
\end{equation*}
In particular, the divisibility of $K_X$ is the greatest common divisor of $n-2$ and $m-2$.
\end{rem}

\begin{rem} Catanese and Wajnryb \cite{Cat2, Cat3, CW} have constructed certain families of simply-connected surfaces of general type with divisible canonical class, using branched coverings over singular curves. Some of these surfaces are diffeomorphic but not deformation equivalent, thus giving counter-examples to a well-known conjecture. 
\end{rem}

\bibliographystyle{amsplain}

\bigskip
\bigskip

\end{document}